 \documentclass[a4paper,11pt]{article}
  \usepackage{amssymb,amsbsy,amsmath,latexsym,theorem,epic,eepic,enumerate}
  \usepackage{xypic,times,overpic}
  \xyoption{curve}
  \usepackage[dvips]{geometry}
  \usepackage{charter,eulervm}
\usepackage[ps2pdf, bookmarks=true,
            bookmarksnumbered=true, breaklinks=true,
            pdfstartview=FitH, hyperfigures=false,
            plainpages=false, naturalnames=true,
            colorlinks=true,pagebackref=true,
            pdfpagelabels]{hyperref}
\usepackage[usenames]{color}

{\theorembodyfont{\rmfamily}\newtheorem{example}{Example}[section]}
{\theorembodyfont{\rmfamily}\newtheorem{Def}[example]{Definition}}
\newtheorem{prop}[example]{Proposition}
\newtheorem{thm}[example]{Theorem}
{\theorembodyfont{\rmfamily}}
\newtheorem{cor}[example]{Corollary}

\newtheorem{blank}[example]{\hspace{-0.45em}}
\newtheorem{diag}[example]{Diagram}

  \def\FTop{\mathsf{FTop}}
 \def\FTOP{\mathsf{FTOP}}
  
  \def\Crs{\mathsf{Crs}}

\def\Set{\mathsf{Set}}
\def\simp{\mathsf{Simp}}
\def\Chain{\mathsf{Chain}}
\def\cU{\mathcal{U}}
\def\geq{\geqslant}
\def\leq{\leqslant}
\def\ecoll{\searrow \hspace{-0.75em}^e\hspace{0.75em}}
\newcommand{\labto}[1]{\stackrel{#1}{\longrightarrow}}

\def\Ce{\mbox{\v{C}}}

\def\Coker{\mathrm{Coker}}
\def\ab{\mathrm{ab}}

\newenvironment{proof}{\noindent {\bf Proof }}{\rule{0mm}{1mm}\hfill $\Box$

\mbox{}}
\newcommand{\sqdiagram}[8]{ \diagram  #1  \rto^-{#2} \dto_{#4} &
#3  \dto^{#5} \\ #6    \rto_-{#7}  &  #8   \enddiagram }
\newcommand{\threeaxes}[3]{\def\objectstyle{\scriptstyle}  \objectmargin={0pt}
\xy
(0,0)*+{}="a",(0,-6)*+{\rule{0em}{1.5ex}#2}="b",(7,0)*+{\;#1}="c",
(14,-3)*+{\;#3}="d" \ar@{->} "a";"b" \ar @{->}"a";"c"  \ar
@{->}"a";"d"\endxy }

\newcommand{\directs}[2]{\def\objectstyle{\scriptstyle}  \objectmargin={0pt}
\xy
(0,4)*+{}="a",(0,-2)*+{\rule{0em}{1.5ex}#2}="b",(7,4)*+{\;#1}="c"
\ar@{->} "a";"b" \ar @{->}"a";"c" \endxy }
\newcommand{\xdirects}[2]{\def\objectstyle{\scriptstyle}  \objectmargin={0pt}
\xy
(0,0)*+{}="a",(0,-6)*+{\rule{0em}{1.5ex}#2}="b",(7,0)*+{\;#1}="c"
\ar@{->} "a";"b" \ar @{->}"a";"c" \endxy }
\newcommand{\sdirects}[2]{\def\objectstyle{\scriptstyle}  \objectmargin={0pt}
\xy
(0,2.2)*+{}="a",(0,-2.5)*+{\rule{0em}{1.5ex}#2}="b",(7,2.2)*+{\;#1}="c"
\ar@{->} "a";"b" \ar @{->}"a";"c" \endxy }

\newcommand{\bl}{\mbox{\rule{0.08em}{1.7ex}\hspace{-0.00em}\rule{0.7em}{0.2ex}}}

\newcommand{\br}{\mbox{\rule{0.7em}{0.2ex}\hspace{-0.04em}\rule{0.08em}{1.7ex}}}

\newcommand{\tr}{\mbox{\rule[1.5ex]{0.7em}{0.2ex}\hspace{-0.03em}\rule{0.08em}{1.7ex}}}

\newcommand{\tl}{\mbox{\rule{0.08em}{1.7ex}\rule[1.54ex]{0.7em}{0.2ex}}}

\newcommand{\hh}{\mbox{\rule{0.7em}{0.2ex}\hspace{-0.7em}\rule[1.5ex]{0.70em}{0.2ex}}}

\newcommand{\vv}{\mbox{\rule{0.08em}{1.7ex}\hspace{0.6em}\rule{0.08em}{1.7ex}}}

\newcommand{\tsq}{\mbox{\rule{0.04em}{1.55ex}\hspace{-0.00em}\rule{0.7em}{0.1ex}\hspace{-0.7em}\rule[1.445ex]{0.7em}{0.1ex}\hspace{-0.03em}\rule{0.04em}{1.55ex}}}

\def\rho{\varrho}

\def\A{\alpha}

\def\I{\mathsf{I}}
\def\eps{\varepsilon}
\def\epsilon{\varepsilon}
\def\pt{\partial}

\def\Z{\mathbb{Z}}

\def\cal{\mathcal}
\def\le{\leqslant}
\def\ge{\geqslant}

\def\x{\mathbf{x}}
\def\y{\mathbf{y}}

\def\d{\mathbf{d}}

\def\lan{\langle}
\def\ran{\rangle}

\def\A{\alpha}

\def\eps{\varepsilon}
\def\epsilon{\varepsilon}
\def\pt{\partial}

\def\tilde{\widetilde}
\def\Z{\mathbb{Z}}

\def\cal{\mathcal}
\def\le{\leqslant}
\def\ge{\geqslant}
\def\leq{\leqslant}
\def\geq{\geqslant}
\def\subset{\subseteq}
\def\ecoll{\searrow \hspace{-0.75em}^e\hspace{0.75em}}

 \def\c{\mathbin{\#}}
 
 \def\d{\partial}

 \def\Im{\mathop{\rm Im}\nolimits}
 
\def\Ker{\mathop{\rm Ker}\nolimits}

 \def\eps{\varepsilon}

\def\crs{\mathsf{Crs}}
\def\Crs{\mathsf{CRS}}

\def\A{\alpha}
\def\bI{\mathbb{I}}
\def\bC{\mathbb{C}}

\def\bZ{\mathbb{Z}}
\def\om{$\omega$-}

\def\ogpd{$\omega$-$\mathsf{Gpd}$}

\def\ogpdm{\omega\mathrm{-}\mathsf{Gpd}}
\def\xybiglabels{\def\labelstyle{\textstyle}}

\def\cal{\mathcal}
\def\epsilon{\varepsilon}

\def\cP{\mathcal{P}}

\newcommand{\bu}{\mathbf{U}}
\def\C{\mathsf{C}}
\newcommand\quadr[4]{\left(\begin{smallmatrix} & #1 &\\ #2 && #3\\&#4& \end{smallmatrix} \right)}
\newcommand{\io}{^{-1}}

\def\blue{\textcolor{blue}}
\def\bu{\bullet}
\begin{document}

\title{Crossed complexes and  higher homotopy groupoids \\ as  non commutative tools for
higher dimensional \\ local-to-global problems\thanks{This is a
revised version (2007) of  a paper published in  Fields Institute
Communications 43 (2004) 101-130, which was an extended account of a
lecture given at the meeting on `Categorical Structures for Descent,
Galois Theory, Hopf algebras and semiabelian categories', Fields
Institute, September 23-28, 2002. The author is grateful for support
from the Fields Institute and a Leverhulme Emeritus Research
Fellowship, 2002-2004, and to M. Hazewinkel for helpful comments on
a draft. }}
\author{Ronald Brown\thanks{ School of
Computer Science, Bangor University, Dean St., Bangor,  Gwynedd LL57
1UT, U.K. email: r.brown{@}bangor.ac.uk}} \maketitle

{\small  MATHEMATICS SUBJECT CLASSIFICATION:
01-01,16E05,18D05,18D35,55P15,55Q05}

\begin{abstract}
We outline the main features of the definitions and applications of
crossed complexes and cubical $\omega$-groupoids with connections.
These give forms of higher homotopy groupoids, and new views of
basic algebraic topology and the cohomology of groups, with the
ability to obtain some non commutative results and compute some
homotopy types in non simply connected situations.
\end{abstract}

 \addcontentsline{toc}{section}{Introduction}
 \tableofcontents

\section*{Introduction}\label{sec:intro}
An aim  is to give a survey and explain the origins  of results
obtained by R. Brown and P.J. Higgins and others over the years
1974-2008, and to point to applications and related areas. These
results  yield an account of some basic algebraic topology on the
border between homology and homotopy; it differs from the standard
account through the use of {\it crossed complexes}, rather than
chain complexes, as a fundamental notion. In this way one obtains
comparatively quickly\footnote{This comparison is based on the fact
that the methods do not require singular homology or simplicial
approximation. } not only classical results such as the Brouwer
degree and the relative Hurewicz theorem, but also non commutative
results on second relative homotopy groups, as well as higher
dimensional results involving the fundamental group, through its
actions  and presentations.  A basic tool is  the fundamental
crossed complex $\Pi X_*$ of the filtered space $X_*$, which in the
case $X_0$ is a singleton is fairly classical; applied to the
skeletal filtration of a $CW$-complex $X$,  $\Pi$ gives a more
powerful version of the usual cellular chains of the universal cover
of $X$, because it contains non-Abelian  information in dimensions 1
and 2, and has good realisation properties. It also gives a
replacement for singular chains by taking $X$ to be the geometric
realisation of a singular complex of a space.

One of the major results is a homotopy classification theorem
(4.1.9) which generalises a classical theorem of Eilenberg-Mac~Lane,
though this does require results on geometric realisations of
cubical sets.

A replacement for the excision theorem in homology is obtained by
using cubical methods to prove a Higher Homotopy van Kampen Theorem
(HHvKT)\footnote{We originally called this a generalised van Kampen
Theorem, but this new term was suggested in 2007 by Jim Stasheff.}
for the fundamental crossed complex functor $\Pi$ on filtered
spaces. This theorem is a higher dimensional version of the van
Kampen Theorem (vKT) on the fundamental group of a space with base
point, \cite{kampen1}\footnote{An earlier version for simplicial
complexes is due to Seifert.}, which is a classical example of a
$$\text{\it non commutative local-to-global theorem},$$ and  was the initial
motivation for the work described here. The vKT determines
completely the fundamental group $\pi_1(X,x)$ of a space $X$ with
base point which is the union of open sets $U,V$ whose intersection
is path connected and contains the base point $x$; the `local
information' is on the morphisms of fundamental groups induced by
the inclusions $U \cap V \to U, U \cap V \to V$. The importance of
this result reflects the importance of the fundamental group in
algebraic topology, algebraic geometry, complex analysis, and many
other subjects. Indeed the origin of the fundamental group was in
Poincar\'e's work on monodromy for complex variable theory.

Essential  to this  use of crossed complexes, particularly for
conjecturing and proving local-to-global theorems,  is a
construction of a {\it cubical higher homotopy groupoid}, with
properties described by {\it an algebra of cubes}. There are
applications to local-to-global problems in homotopy theory which
are more powerful than available by purely classical tools, while
shedding light on those tools. It is hoped that this account will
increase the interest in the possibility of wider applications of
these methods and results, since homotopical methods play a key role
in many areas.

\noindent {\bf Background in higher homotopy groups}

\noindent  Topologists in the early part of the 20th century were
well aware that: \begin{itemize}
\item the non commutativity of the
fundamental group was useful in geometric applications;
\item for path
connected $X$ there was an isomorphism
$$H_1(X) \cong \pi_1(X,x)^{\mathrm{ab}}; $$
\item the Abelian  homology groups $H_n(X)$ existed for all $n \geq 0$.
\end{itemize}
Consequently there was a desire to generalise the non commutative
fundamental group to all dimensions.

In 1932 \v{C}ech submitted a paper on higher homotopy groups
$\pi_n(X,x)$ to the ICM at Zurich, but it was quickly proved that
these groups were Abelian  for $n \ge 2$, and on these grounds
\v{C}ech was persuaded to withdraw his paper, so that only a small
paragraph appeared in the Proceedings \cite{Cech1}. We now see the
reason for this commutativity as the result (Eckmann-Hilton) that a
group internal to  the category of groups is just an Abelian  group.
Thus, since 1932 the vision of a non commutative higher dimensional
version of the fundamental group has been generally considered to be
a mirage. Before we go back to the vKT, we explain in the next
section  how nevertheless work on crossed modules did introduce non
commutative structures relevant to topology in dimension 2.

Work of Hurewicz, \cite{hurewicz}, led to a strong development of
higher homotopy groups. The fundamental group still came into the
picture with its action on the higher homotopy groups, which I once
heard J.H.C. Whitehead remark (1957) was especially fascinating for
the early workers in homotopy theory. Much of Whitehead's work was
intended to extend to higher dimensions the methods of combinatorial
group theory of the 1930s -- hence the title of his papers:
`Combinatorial homotopy, I, II' \cite{jhcw:CHI,jhcw:CHII}. The first
of these two papers has been very influential and is part of the
basic structure of algebraic topology. It is the development of work
of the second paper which we explain here.

The paper by Whitehead on `Simple homotopy types' \cite{wjhc:sht},
which deals with higher dimensional analogues of Tietze
transformations, has a final section using crossed complexes. We
refer to this again later in section \ref{S:freecrossed}.

It is hoped also that this survey will be useful background to work
on the van Kampen Theorem for diagrams of spaces in \cite{brlo:vkt},
which uses a form of higher homotopy groupoid which is in an
important  sense much more powerful than that given here, since it
encompasses $n$-adic information; however current expositions are
still restricted to the reduced (one base point) case, the proof
uses advanced tools of algebraic topology, and the result was
suggested by the work exposed here.

\section{Crossed modules}
\label{sec:crossedmod} In the years 1941-50, Whitehead  developed
work on crossed modules to represent the structure of the boundary
map of the relative homotopy group
\begin{equation} \pi_2(X,A,x) \to \pi_1(A,x) \label{relhom} \end{equation}
in which both groups can be non commutative. Here is the definition
he found.

A {\em crossed module} is a morphism of groups $\mu : M \to P$
together with an action $(m,p)\mapsto m^p$ of the group $P$ on the
group $M$ satisfying the two axioms
\begin{enumerate}[{CM}1)] \item $\mu
(m^p) = p^{-1}(\mu m)p$ \item $n^{-1}mn= m ^{\mu n}$
\end{enumerate} for all $m,n \in M, p \in P.$
\par Standard algebraic examples of crossed modules are:
\begin{enumerate}[(i)]
\item an inclusion of a normal subgroup, with action given by
conjugation;
\item the inner automorphism map $\chi :
M \to \mbox{Aut}\;M,$ in which $\chi m $ is the automorphism
$n\mapsto m^{-1}nm$; \item the zero  map $M \to P $ where  $M$ is a
$P$-module; \item an epimorphism $M \to P $ with kernel contained in
the centre of $M$.
\end{enumerate}
Simple consequences of the axioms for a crossed module $\mu : M \to
P$ are:
\begin{blank} $\Im \mu $ is normal in $P$. \end{blank}
\begin{blank} $\Ker \mu$ is central in $M$ and is acted on trivially
by $\Im \mu$, so that $\Ker \mu$ inherits an action of $M/\Im
\mu$.
\end{blank}
Another important algebraic construction is the {\em free crossed
$P$-module }
$$\partial : C(\omega) \to P $$ determined by a function $\omega: R \to P$,
where $P$ is a group and $R$ is a set. The group  $C(\omega)$ is
generated by elements $(r,p) \in R\times P$ with the relations
$$(r,p)\io (s,q)\io (r,p) (s,qp\io(\omega r)p);$$ the action is
given by $(r,p)^q=(r,pq)$; and the boundary morphism is given by
$\partial (r,p) = p\io(\omega r)p$, for all $(r,p),(s,q) \in R\times
P$.

A major result of Whitehead was:

\noindent {\bf Theorem W}  \cite{jhcw:CHII} {\em If the space $X=A
\cup \{e^2_r\}_{r \in R}$ is obtained from $A$ by attaching
2-cells by maps $f_r:(S^1,1) \to (A,x)$, then the crossed module
of \eqref{relhom} is isomorphic to the free crossed
$\pi_1(A,x)$-module on the classes of the attaching maps of the
2-cells.}

Whitehead's proof, which stretched over three papers, 1941-1949,
used transversality and knot theory -- an exposition is given in
\cite{rb:second80}. Mac~Lane and Whitehead \cite{maclane&jhcw}
used this result as part of their proof that crossed modules
capture all homotopy 2-types (they used the term `3-types').

The title of the paper in which the first intimation of Theorem W
appeared was `On adding relations to homotopy groups'
\cite{whjhc:adding41}. This indicates a search for higher
dimensional vKTs.

The concept of free crossed module gives a non commutative context
for {\it chains of syzygies}. The latter idea, in the case of
modules over polynomial rings, is one of the origins of
homological algebra through the notion of {\it free resolution}.
Here is how similar ideas can be applied to groups. Pioneering
work here, independent of Whitehead, was by Peiffer \cite{peiffer}
and Reidemeister \cite{reidemeister}. See \cite{brownhuebschmann}
for an exposition  of these ideas.

Suppose $\cP= \lan X \mid\omega\ran $ is a presentation of a group
$G$, so that $X$ is a set of generators of $G$ and  $\omega: R \to
F(X) $ is a function, whose image is called the set of relators of
the presentation. Then we have an exact sequence
$$1 \labto{i} N(\omega R) \labto{\phi} F(X) \longrightarrow G \longrightarrow 1$$ where $N(\omega R)$ is
the normal closure in $F(X)$ of the set $\omega R$ of relators. The
above work of Reidemeister, Peiffer, and Whitehead showed that to
obtain the next level of syzygies one should consider the free
crossed $F(X)$-module $\partial: C(\omega) \to F(X)$, since this
takes into account  the operations of $F(X)$ on its normal subgroup
$N(\omega R)$. Elements of $C(\omega)$ are a kind of `formal
consequences of the relators', so that the relation between the
elements of $C(\omega)$ and those of $N(\omega R)$ is analogous to
the relation between the elements of $F(X)$ and those of $G$.
\blue{It follows from the rules for a crossed module that the kernel
of $\partial$ is a $G$-module, called the module of {\it identities
among relations}, and sometimes written $\pi(\cP)$; there is
considerable work on computing it
\cite{brownhuebschmann,pride,hog-ang,ellis-kholod,BRazak:LMS99}}. By
splicing to $\partial$ a free $G$-module resolution of $\pi(\cP)$
one obtains what is called a {\it free crossed resolution } of the
group $G$. We explain later (Proposition \ref{cwmaps}) why these
resolutions have better realisation properties than the usual
resolutions by chain complexes of $G$-modules. They are relevant to
the Schreier extension theory, \cite{brposchreier}.

This notion of using crossed modules as the first stage of
syzygies in fact represents a wider tradition in homological
algebra, in the work of Fr\"{o}lich and Lue
\cite{frohlich:nonab,lue2}.

Crossed modules also occurred in other contexts, notably in
representing elements of the cohomology group $H^3(G,M)$ of a group
$G$ with coefficients in  a $G$-module $M$ \cite{maclane:hom}, and
as coefficients in Dedecker's theory of non Abelian  cohomology
\cite{ded-canj}. The notion of  free crossed resolution has been
exploited by Huebschmann
\cite{huebschmann,huebschmann2,huebschmann3} to represent cohomology
classes in $H^n(G,M)$ of a group $G$ with coefficients in a
$G$-module $M$, and also to calculate with these.

The HHvKT can make it easier to compute a crossed module arising
from some topological situation, such as an induced crossed module
\cite{brownwens:ind,brownwens:normal}, or a coproduct crossed module
\cite{rb:coproducts}, than the cohomology class in $H^3(G,M)$
\blue{the crossed module}  represents. To obtain information on such
a cohomology element  it is useful to work with a small free crossed
resolution of $G$, and this is one motivation for developing methods
for calculating such resolutions. However, it is not so clear what a
{\it calculation} of such a cohomology element would amount to,
although it is interesting to know whether the element is non zero,
or what is its order. Thus the use of algebraic models of cohomology
classes may yield easier computations than the use of cocycles, and
this somewhat inverts traditional approaches.

Since crossed modules are algebraic objects generalising groups,
it is natural to consider the problem of explicit calculations by
extending  techniques of computational group theory. Substantial
work on this has been done by C.D. Wensley using the program GAP
\cite{GAP,brownwens:comp}.

\section{The fundamental groupoid on a set of base points}
\label{sec:fundgroupoid} A change in prospects for higher order non
commutative invariants was suggested by Higgins' paper
\cite{Higgins2}, and  leading  to work of the writer published in
1967, \cite{brown:gpdsvkt67}. This showed that the van Kampen
Theorem could be formulated for the {\it fundamental groupoid
$\pi_1(X,X_0)$ on a {\it set} $X_0$ of base points}, thus enabling
computations in the non-connected case, including those in Van
Kampen's original paper \cite{kampen1}. This successful use of
groupoids in dimension 1 suggested the question of the use of
groupoids in higher homotopy theory, and in particular the question
of the existence of {\it higher homotopy groupoids}.

In order to see how this research programme could progress it is
useful to consider the statement and special features of this
generalised van Kampen Theorem for the fundamental groupoid. If
$X_0$ is a set, and $X$ is a space, then $\pi_1(X,X_0)$ denotes the
fundamental groupoid on the set $X \cap X_0$ of base points. This
allows the set $X_0$ to be chosen in a way appropriate to the
geometry. For example, if the circle $S^1$ is written as the union
of two semicircles $E_+ \cup E_{-}$, then the intersection
$\{-1,1\}$ of the semicircles is not connected, so it is not clear
where to take the base point. Instead one takes $X_0=\{-1,1\}$, and
so has two base points. This flexibility is very important in
computations, and this example of $S^1$ was a motivating example for
this development. As another example, you might like to consider the
difference between the quotients of the actions of $\bZ_2$ on the
group $\pi_1(S^1,1)$ and on the groupoid $\pi_1(S^1,\{-1,1\})$ where
the action is induced by complex conjugation on $S^1$. Relevant work
on {\it orbit groupoids} has been developed by Higgins and Taylor
\cite{HigginsTaylor1,taylorj:88}, (under useful conditions, the
fundamental groupoid of the orbit space is the orbit groupoid of the
fundamental groupoid \cite[11.2.3]{brownbook:2}).

Consideration of a set of base points leads to the theorem:
\begin{thm}{\em \cite{brown:gpdsvkt67}}\label{vktgpd} Let the space $X$ be the
union of open sets $U,V$ with intersection $W$, and let $X_0$ be a
subset of $X$
meeting each path component of $U,V,W$. Then \\
 {\em (C) (connectivity)}  $X_0$ meets each path component of $X$ and \\
 {\em (I) (isomorphism)}  the diagram of groupoid morphisms induced by inclusions
\begin{equation} \label{push}
{\sqdiagram{\pi_1(W,X_0)}{i}{\pi_1(U,X_0)}{j}{{j'}}{\pi_1(V,X_0)}{{i'}}{\pi_1(X,X_0)}}
\end{equation}
is a pushout of groupoids. \end{thm}

From this theorem, one can compute a particular fundamental group
$\pi_1(X,x_0)$ using combinatorial information on the graph of
intersections of path components of $U,V,W$, but for this it is
useful  to develop the algebra of groupoids. Notice two special
features of this result.

\noindent (i) The computation of the invariant you may want, a
fundamental group, is obtained from the computation of a larger
structure, and so part of the work is to give methods for computing
the smaller structure from the larger one. This usually involves non
canonical choices, e.g. that of a maximal tree in a connected graph.
The work on applying groupoids to groups gives many examples of this
\cite{Higgins2,Higgins4,brownbook:2,dicks-vent}.

\noindent (ii) The fact that the computation can be done is
surprising in two ways: (a) The fundamental group is computed {\it
precisely}, even though the information for it uses input in two
dimensions, namely 0 and 1. This is contrary to the experience in
homological algebra and algebraic topology, where the interaction
of several dimensions involves exact sequences or spectral
sequences, which give information only up to extension,  and (b)
the result is a non commutative invariant, which is usually even
more difficult to compute precisely.

The reason for the success seems to be that the fundamental groupoid
$\pi_1(X,X_0)$ contains information in dimensions 0 and 1, and so
can adequately reflect the geometry of the intersections of the path
components of $U,V,W$ and of the morphisms induced by the inclusions
of $W$ in $U$ and $V$.

This suggested the question of  whether these methods could be
extended successfully to higher dimensions.

Part of the initial evidence for this quest was the intuitions in
the proof of this groupoid vKT, which seemed to use  three main
ideas in order to verify the universal property of a pushout for
diagram \eqref{push}. So suppose given morphisms of groupoids
$f_U,f_V$ from $\pi_1(U,X_0), \pi_1(V,X_0)$ to a groupoid $G$,
satisfying $f_Ui=f_Vj$. We have to construct a morphism
$f:\pi_1(X,X_0) \to G$ such that $fi'=f_U,fj'=f_V$ and prove $f$ is
unique. We concentrate on the construction.

$\bullet$ One needs a `deformation', or `filling', argument: given a
path $a:(I, \dot I) \to (X,X_0)$ one can write $a=a_1+\cdots +a_n$
where each $a_i$ maps into $U$ or $V$, but $a_i$ will not
necessarily have end points in $X_0$. So one has to deform each
$a_i$ to $a'_i$ in $U,V$ or $W$, using the connectivity condition,
so that each $a'_i$ has end points in $X_0$, and $a'=a'_1+\cdots
+a'_n$ is well defined. Then one can construct using $f_U$ or $f_V$
an image of each $a'_i$ in $G$ and hence of the composite, called
$F(a) \in G$, of these images. Note that we subdivide in $X$ and
then put together again in $G$ (this uses the condition $f_Ui=f_Vj$
to prove that the elements of $G$ are composable), and this part can
be summarised as:

$\bullet$ Groupoids provided a convenient {\em algebraic inverse to
subdivision}. Note that the usual exposition in terms only of the
fundamental group uses loops, i.e. paths which start and finish at
the same point. An appropriate analogy  is that  if one goes on a
train journey from Bangor and back to Bangor, one usually wants to
stop off at intermediate stations; this breaking and cotinuing a
journey is better described in terms of groupoids rather than
groups.

Next one has to prove that $F(a)$ depends only on the class of $a$
in the fundamental groupoid. This involves a homotopy rel end points
$h:a \simeq b$, considered as a map $I^2 \to X$; subdivide $h$ as
$h=[h_{ij}]$ so that each $h_{ij}$ maps into $U,V$ or $W$; deform
$h$ to $h'=[h'_{ij}]$ (keeping in $U,V,W$) so that each $h'_{ij}$
maps the vertices to $X_0$ and so determines a commutative
square\footnote{\blue{We need the notion of  {\it commutative
square} in a category $\,\C$. This  is a quadruple
$\quadr{c}{a}{d}{b}$ of arrows in $\C$, called `edges' of the
square, such that $ab=cd$, i.e. such that these compositions are
defined and agree. The commutative squares in $\C$ form a {\it
double category} $\square \C$ in that they compose `vertically'
$$\quadr{c}{a}{d}{b} \circ_1 \quadr{b}{a'}{d'}{e} =
\quadr{c}{aa'}{dd'}{e}$$and `horizontally'
$$\quadr{c}{a}{d}{b} \circ_2 \quadr{c'}{d}{f}{b'} =
\quadr{cc'}{a}{f}{bb'}$$ This notion of $\square \C$ was defined by
C. Ehresmann in papers and in \cite{Eh-structure}. Note the obvious
geometric conditions for these compositions to be defined.
Similarly, one has geometric conditions for  a rectangular array
$(c_{ij}), 1 \leq i \leq m,1 \leq j \leq n,  $ of commutative
squares to  have a well defined composition, and then their
`multiple composition', written $[c_{ij}]$, is also a commutative
square, whose edges are compositions of the `edges' along the
outside boundary of the array. It is easy to give formal definitions
of all this.  }} in one of $\pi_1(Q,X_0)$ for $Q=U,V,W$. Move these
commutative squares over to $G$ using $f_U,f_V$ and recompose them
(this is possible again because of the condition $f_Ui=f_Vj$),
noting that:

$\bullet$ in a groupoid, {\em any  composition of commutative
squares is commutative}. Here a `big' composition of commutative
squares is represented by a diagram such as
\begin{equation}\label{eq:multcomp}{\objectmargin{0.1pc} \diagram \bu \rto \dto &
\bu \rto \dto & \bu \rto \dto & \bu \rto \dto & \bu \rto \dto & \bu
\rto \dto & \bu  \dto \\    \bu \rto \dto &  \bu \rto \dto & \bu
\rto \dto & \bu \rto \dto & \bu \rto \dto & \bu \rto \dto & \bu \dto
\\  \bu \rto \dto &  \bu \rto \dto  & \bu \rto \dto & \bu \rto \dto & \bu \dto\rto & \bu
\rto \dto & \bu  \dto \\ \bu \rto \uto  & \bu \rto \uto
 & \bu \rto    & \bu \rto    & \bu \rto   \ & \bu
\rto   & \bu    \enddiagram  }\end{equation} and one checks that if
each individual square is commutative, so also is the boundary
square (later called a 2-shell) of the compositions of the boundary
edges.

\noindent Two opposite  sides of the composite commutative square in
$G$ so obtained are identities, because $h$ was a homotopy relative
to end points, and the other two sides are $F(a),\, F(b)$.  This
proves that $F(a)=F(b)$ in $G$.

Thus the argument can be summarised:  a path or homotopy is
divided into small pieces, then  deformed so that these pieces can
be packaged and moved over to $G$, where they are reassembled.
There seems to be  an analogy with the processing of an email.

Notable applications of the groupoid theorem were: (i) to give a
proof of a formula in van Kampen's paper of the fundamental group of
a space which is the union of two connected spaces with non
connected intersection,  see \cite[8.4.9]{brownbook:2}; and (ii) to
show the topological utility of the construction by Higgins
\cite{Higgins4} of the groupoid $f_*(G)$ over $Y_0$ induced from a
groupoid $G$ over $X_0$ by a function $f:X_0 \to Y_0$. (Accounts of
these with the notation $U_f(G)$ rather than $f_*(G)$ are given in
\cite{Higgins4,brownbook:2}.) This latter construction is regarded
as a `change of base', and analogues  in higher dimensions yielded
generalisations of the Relative Hurewicz Theorem and of Theorem W,
using induced modules and crossed modules.

There is another approach to the van Kampen Theorem which goes via
the theory of covering spaces, and the equivalence between covering
spaces of a reasonable space $X$ and functors $\pi_1(X) \to \Set$
\cite{brownbook:2}. See for example \cite{Douady1b} for an
exposition of the relation with traditional Galois theory, and
\cite{borc-janelid} for a modern account in which {\it Galois
groupoids} make an essential appearance. The paper
\cite{brownjan:vkt} gives a general formulation of conditions for
the theorem to hold in the case $X_0=X$ in terms of the map $U
\sqcup V \to X$ being an `effective global descent morphism' (the
theorem is given in the generality of lextensive categories). This
work has been developed for toposes, \cite{BungeL}. Analogous
interpretations for higher dimensional Van Kampen theorems are not
known.

The justification of the breaking of a paradigm in changing from
groups to groupoids is several fold: the elegance and power of the
results; the increased linking with other uses of groupoids
\cite{brown:survey}; and the opening out of new possibilities in
higher dimensions, which allowed for new results and calculations in
homotopy theory, and suggested new algebraic constructions. The
important and extensive work of Charles Ehresmann in using groupoids
in geometric situations (bundles, foliations, germes, $\ldots$)
should also be stated (see his collected works of which
\cite{ehres-works} is volume 1 and a survey \cite{brown:ehres}).

\section{The search for higher homotopy groupoids}
\label{sec:search} Contemplation of  the proof of the groupoid vKT
in the last section suggested that a higher dimensional version
should exist, though this version amounted to an idea of a proof in
search of a theorem. Further evidence was the proof by J.F. Adams of
the cellular approximation theorem given in \cite{brownbook:2}. This
type of subdivision argument failed to give algebraic information
apparently because of a lack  of an appropriate higher homotopy
groupoid, i.e. a gadget to capture what might be the underlying
`algebra of cubes'. In the end, the results exactly encapsulated
this intuition.

One intuition was that in groupoids we are dealing with a partial
algebraic structure\footnote{The study of partial  algebraic
operations was initiated in \cite{Higgins1}. We can  now suggest a
reasonable definition of `higher dimensional algebra' as dealing
with families of algebraic  operations whose domains of definitions
are given by geometric conditions.}, in which composition is defined
for two arrows if and only if the source of one arrow is the target
of the other. This seems  to generalise easily to directed squares,
in which two such are composable horizontally if and only if the
left hand side of one is the right hand side of the other (and
similarly vertically).

However the formulation of a theorem in higher dimensions required
specification of  the {\it topological data}, the  {\it algebraic
data}, and of a functor
$$\Pi: \mbox{(topological data)} \to \mbox{(algebraic data)}$$
which would allow the expression of these ideas for the
proof.

Experiments were made in the years 1967-1973 to define some functor
$\Pi$ from spaces to some kind of double groupoid, using
compositions of squares in two directions, but these proved
abortive. However considerable progress was made in work with Chris
Spencer in 1971-3 on investigating the algebra of double groupoids
\cite{brownspenc:double}, and showing a relation to crossed modules.
Further evidence was provided when it was found,
\cite{brownspenc:g-gr}, that group objects in the category of
groupoids (or groupoid objects in the category of groups, either of
which are often now called `2-groups') are equivalent to crossed
modules, and in particular are not necessarily commutative objects.
It turned out this result was known to the Grothendieck school in
the 1960s, but not published.

{We review next a notion of double category which is not the most
general but is appropriate in  many cases.}  It  was  called an {\it
edge symmetric double category} in \cite{brownmosa}.

In the first place, a double category, $ { \mathsf K} $, consists of
a triple of category structures
\[
\begin{array}{c}
(K_{2} ,K_{1} ,\partial^{-}_{1} , \partial^{+}_{1} , \circ_{1}
,\varepsilon_{1} ) , \quad
(K_{2} , K_{1} ,\partial^{-}_{2} , \partial^{+}_{2} ,\circ _{2} ,\varepsilon_{2} )\\[2mm]
(K_{1} , K_{0} ,\partial^{-} , \partial^{+} , \circ  , \varepsilon )
\end{array}
\]
 as partly shown in the diagram
\begin{equation}\xybiglabels
\xymatrix@=4pc{K_2 \ar @<0.5ex> [d]^{\partial^-_1} \ar @<-0.5ex>
[d]_{\partial^+_1} \ar @<0.5ex> [r]^{\partial^-_2} \ar @<-0.5ex>
[r]_{\partial^+_2}& K_1  \ar @<0.5ex> [d]^{\partial^-} \ar
@<-0.5ex> [d]_{\partial^+} \\
K_1 \ar @<0.5ex> [r]^{\partial^-} \ar @<-0.5ex> [r]_{\partial^+}&
K_0}
\end{equation}
The elements of $K_0,K_1,K_2$ will be called respectively {\it
points or objects, edges, squares}. The maps $ {\partial^\pm,} \;
\partial^{\pm}_{i} , \ i = 1,2 $, will be called {\it face maps},
the maps $ \varepsilon_{i} : K_{1} \longrightarrow K_{2} , \ i = 1,2
$, resp. $ \varepsilon : K_{0} \longrightarrow K_{1} $ will be
called {\it degeneracies}. The boundaries of an edge and of a square
are given by the diagrams
\begin{equation}
\vcenter{\xymatrix@=3pc{\partial^- \ar @{-}@<-0.3ex>[r]|\tip &
\partial^+}}\qquad
 \vcenter{\xymatrix@=3pc@M=0pt{ \ar @{-}[r]|\tip ^{\partial^-_1}\ar @{-}[d]_{\partial^-_2}|\tip & \ar @{-}[d]|\tip  ^{\partial^+_2}\\
  \ar @{-}[r]_{\partial^+_1}|\tip  &    }}\quad \directs{2}{1}
\end{equation}
The partial compositions, $ \circ _{1} $, resp. $ \circ _{2} $, are
referred to as {\it vertical} resp. {\it horizontal composition} of
squares, are defined under the obvious geometric conditions, and
have the obvious boundaries. The axioms for a double category also
include the usual relations of a 2-cubical set (for example
$\partial^-\partial^+_2=\partial^+\partial^-_1$), and the {\it
interchange law}. We use matrix notation for compositions as
$$\begin{bmatrix}a\\c \end{bmatrix} = a \circ _1 c, \quad
 \begin{bmatrix}a & b \end{bmatrix} = a \circ _2 b, $$ and the crucial  interchange law\footnote{The interchange
implies that a double monoid is simply an Abelian monoid, so partial
algebraic operations are essential for the higher dimensional work.}
for these two compositions  allows one to use matrix notation
$$\begin{bmatrix} a&b \\c&d \end{bmatrix} = \begin{bmatrix}
\begin{bmatrix}   a&b \end{bmatrix}\\[1ex] \begin{bmatrix}
  c&d
\end{bmatrix}
\end{bmatrix} = \begin{bmatrix}\begin{bmatrix}
  a\\c
\end{bmatrix}& \begin{bmatrix}
  b\\d
\end{bmatrix}
\end{bmatrix}$$ for double composites of squares whenever each row composite and
each column composite is defined. We also allow the multiple
composition $[a_{ij}]$ of an array $(a_{ij})$ whenever for all
appropriate $i,j$ we have $\partial ^+_1a_{ij}=
\partial ^-_1a_{i+1,j},\;\partial  ^+_2a_{ij}= \partial
^-_2a_{i,j+1}$. A clear  advantage of double categories and cubical
methods is this easy expression of multiple compositions which
allows for {\it algebraic inverse to subdivision}, and so
applications to local-to-global problems.

 The identities with respect to $ \circ _{1} $ ({\it vertical
identities}) are given by $ \varepsilon_{1} $ and will be denoted by
$\vv$.  Similarly, we have {\it horizontal identities} denoted by
$\hh$. Elements of the form $ \varepsilon_{1} \varepsilon (a) =
\varepsilon_{2} \varepsilon (a) $ for $ a \in K_{0} $ are called
{\it double degeneracies} and will be denoted by $ \tsq $ .

A \textit{morphism of double categories} $f : \mathsf{K}\to
\mathsf{L}$ consists of a triple of maps $f_i : K_i \to L_i$, $(i =
0,1,2)$, respecting the cubical structure, compositions and
identities.

Whereas it is easy to describe a commutative square of morphisms in
a category, it is not possible with this amount of structure to
describe a commutative cube of squares in a double category. We
first of all define a cube, or 3-shell, i.e. without any condition
of commutativity, in a double category.

\begin{Def} Let $ { \mathsf K} $ be a double
category. A {\it cube}  { (\emph{$3$-shell})} in $ { \mathsf K} $,
\[
\alpha = (\alpha^{-}_{1} , \alpha^{+}_{1} , \alpha^{-}_{2} ,
\alpha^{+}_{2} , \alpha^{-}_{3} , \alpha^{+}_{3} )
\]
 consists of squares $ \alpha^{\pm}_{i} \in K_{2} \quad (i = 1,2,3) $ such that
\[
\partial^{\sigma}_{i} ( \alpha^{\tau}_{j} ) = \partial^{\tau}_{j-1} (\alpha^{\sigma}_{i} )
\]
 for $ \sigma , \tau = \pm  $ and $ 1 \leq i < j \leq 3 $. \hfill
 $\blacksquare$
\end{Def}

It is also convenient to have the corresponding notion of square, or
2-shell, of arrows in a category. The obvious compositions also
makes these into a double category.

It is not hard to define three compositions of cubes in a double
category so that these cubes form a triple category: this is done in
\cite{bkp-vKT}, or more generally in Section 5 of \cite{bh:algcub}.
A key point is that to define the notion of a {\it commutative cube}
we need extra structure on a double category. Thus this step up a
dimension is non trivial, as was first observed in the groupoid case
in \cite{bh1978}. The problem is that a cube has six faces, which
easily divide into three even and three odd faces. So we cannot  say
as we might like that `the cube is commutative if the composition of
the even faces equals the composition of the odd faces', since there
are no such valid compositions.

The intuitive reason for the need of a new basic structure in that
in a 2-dimensional situation we also need to use the possibility of
`turning an edge clockwise or anticlockwise'. The structure to do
this is as follows.

 A {\it connection pair} on a double category $ { \mathsf K} $
 is given by a pair of maps
\[
\Gamma^{-} , \Gamma^{+} : K_{1} \longrightarrow K_{2}
\]
 whose edges are given by the following diagrams for $ a
\in K_{1} $:

\begin{equation*}
  \Gamma^-(a)= \vcenter{\xybiglabels\xymatrix@M=0pt@=3pc{\ar@{-} [r]^a|(0.6)\tip \ar@{-} [d]_a|(0.6)\tip & \ar@{~}[d]\ar@{-} [d] ^1 \\
  \ar@{~}[r] \ar @{-} [r] _1&}} = \vcenter{\xybiglabels\xymatrix@M=0pt@=3pc{\ar@{}[dr]|\br \ar@{-}[r]|(0.6)\tip  \ar@{-} [d]|(0.6)\tip  & \ar@{~}[d]\ar@{-} [d]  \\
  \ar@{~}[r] \ar @{-} [r] &}}\quad  = \br \qquad  \directs{2}{1}
\end{equation*}

\begin{equation*}
  \Gamma^+(a)= \vcenter{\xybiglabels\xymatrix@M=0pt@=3pc{\ar@{-} [r]^1 \ar@{~}[r] \ar@{-} [d]_1 \ar@{~}[d]&
  \ar@{-} [d] ^a|(0.6)\tip \\
  \ar @{-} [r] _a|(0.6)\tip&}} =
 \raisebox{3.8ex}{\xybiglabels\xymatrix@M=0pt@=3pc{\ar@{}[dr]|\tl  \ar@{~}[r] \ar@{-}[r]|(0.6)\tip  \ar@{-} [d]\ar@{~}[d]  & \ar@{-} [d]^a|(0.6)\tip  \\
   \ar @{-} [r]_a|(0.6)\tip &}}\quad  = \tl \qquad  \directs{2}{1}
\end{equation*}
This `hieroglyphic' notation, which was introduced in
\cite{br-hdgt}, is useful for expressing the laws these connections
satisfy. The first is a pair of cancellation laws which read
\begin{equation*}
  \begin{bmatrix}\, \tl \,  \\ \br  \end{bmatrix} = \hh \; , \qquad
  \begin{bmatrix}\tl\, & \br  \,  \end{bmatrix} = \vv\;,
\end{equation*}
which can be understood as `if you turn right and then left, you
face the same way', and similarly the other way round. They were
introduced in \cite{spencer}. Note that in this matrix notation we
assume that the edges of the connections are such that the
composition is defined.

Two other laws relate the connections to the compositions and read
\begin{equation*}
\begin{bmatrix}\, \tl & \hh\;  \\ \vv & \tl \end{bmatrix} = \tl \, ,
\quad\begin{bmatrix}\, \br & \vv\;  \\ \hh & \br \end{bmatrix} = \br
\; .
\end{equation*}
These can be interpreted as `turning left (or right) with your arm
outstretched is the same as turning left (or right)'. The term
`connections' and the name `transport laws' was because these laws
were suggested by the laws  for path connections in differential
geometry, as explained in \cite{brownspenc:double}. It was proved in
\cite{brownmosa} that a connection pair on a double category $K$ is
equivalent to a `thin structure', namely a morphism of double
categories $\Theta: \square K_1 \to K$ which is the identity on the
edges. The proof requires some `2-dimensional rewriting' using the
connections.

We can now explain what is a `commutative cube' in a double category
$K$ with connection pair.

\begin{Def} \label{hal} Suppose given, { in a double category with
connections $ {\mathsf K} $, a cube (3-shell)}
$$
\alpha = (\alpha^{-}_{1} , \alpha^{+}_{1} , \alpha^{-}_{2} ,
\alpha^{+}_{2} , \alpha^{-}_{3} , \alpha^{+}_{3} ).
$$
 We define the {\it composition of the odd faces} of $\alpha$  to
be
\begin{align} \boldsymbol{\partial}^{\mathrm{odd}} \alpha &=
\begin{bmatrix}\;  \tl & \alpha^-_1 & \br \;\\ \alpha^-_3 & \alpha^+_2
& \hh
\end{bmatrix} \\
\intertext{and the {\it composition of the even faces} of $\alpha$
to be}\boldsymbol{\partial}^{\mathrm{even}} \alpha &=
\begin{bmatrix}\hh & \alpha^-_2 & \alpha^+_3 \\
\tl & \alpha^+_1 & \br\;
\end{bmatrix}
\end{align}
 We define $\alpha$ to be {\it commutative} if it satisfies the
Homotopy Commutativity Lemma (HCL), i.e. \begin{equation}
\boldsymbol{\partial}^{\mathrm{odd}} \alpha=
\boldsymbol{\partial}^{\mathrm{even}} \alpha. \tag{HCL}
\end{equation}
This definition can be regarded as a cubical, categorical (rather
than groupoid) form of the Homotopy Addition Lemma (HAL) in
dimension 3.
\end{Def}

You should draw a 3-shell, label all the edges with letters, and see
that this equation makes sense in that the 2-shells of each side of
equation (HCL) coincide. Notice however that these 2-shells do not
have coincident partitions along the edges: that is the edges of
this 2-shell in direction 1  are formed from different compositions
of the type $ 1 \circ a $ and $a\circ 1$. This definition is
discussed in more detail in \cite{bkp-vKT}, is related to other
equivalent definitions, and it is proved that compositions of
commutative cubes in the three possible directions  are also
commutative. These results are extended to all dimensions in
\cite{hig-thin}; this requires the full structure indicated in
section \ref{sec:omegagroupoid} and also the notion of {\it thin
element} indicated in section \ref{sec:thinelements}.

The initial discovery of connections arose in
\cite{brownspenc:double} from relating crossed modules to double
groupoids. The first example of a double groupoid was the double
groupoid $\tsq G$ of commutative squares in a group $G$. The first
step in  generalising this construction was to  consider quadruples
$\left(\begin{smallmatrix} &c&\\a&&d\\&b&
\end{smallmatrix}\right)$ of elements of $G$ such that $abn=cd$ for some
element $n$ of a subgroup $N$ of $G$. Experiments quickly showed
that for the two compositions of such quadruples to be valid it was
necessary and sufficient that  $N$ be normal in $G$.  But in this
case the element $n$ is determined by the boundary, or 2-shell,
$a,b,c,d$.  In homotopy theory we require something more general. So
we consider a morphism $\mu:N \to G$ of groups and and consider
quintuples $\left( n:\begin{smallmatrix} &c&\\a&&d\\&b&
\end{smallmatrix}\right)$ such that $ab\mu(n)=cd$. It then turns out that we get a
double groupoid if and only if $\mu: N \to G$ is a crossed module.
The next question is which double groupoids arise in this way? It
turns out that we need exactly double groupoids with connection
pairs, though in this groupoid case we can deduce $\Gamma^-$ from
$\Gamma^+$ using inverses in each dimension. This gives the main
result of \cite{brownspenc:double}, the equivalence between the
category of crossed modules and that of double groupoids with
connections and one vertex.

These connections were also used in \cite{bh1978} to define a
`commutative cube' in a double groupoid with connections using the
equation
$$c_1=\begin{bmatrix}\tl & \; a ^{-1}_{0} & \tr \;\\-b_{0} & c_{0} &
b_{1} \\ \bl & a_{1} & \br\end{bmatrix}$$ representing one face of a
cube in terms of the other five and where the other connections
$\bl\, ,\tr $ are obtained from $\br\, , \tl$ by using the two
inverses in dimension 2. As you might imagine, there are problems in
finding a formula  in still higher dimensions. In the groupoid case,
this is handled by a homotopy addition lemma and thin elements,
\cite{bh:algcub}, but in the category case a formula for just a
commutative $4$-cube is complicated, see \cite{gaucher:tac}.

The blockage of defining a functor $\Pi$ to double groupoids was
resolved after 9 years in 1974 in discussions with Higgins, by
considering the Whitehead Theorem W. This showed that a
2-dimensional universal property was available in homotopy theory,
which was encouraging; it also suggested that a theory to be any
good should recover Theorem W. But this theorem was about {\it
relative} homotopy groups. This suggested studying a relative
situation $X_*: X_0 \subseteq X_1 \subseteq X$. On looking for the
simplest way to get a homotopy functor from this situation using
squares,  the `obvious' answer came up: consider maps
$(I^{2},\partial I^2,\partial \partial I^2) \to (X,X_{1},X_{0})$,
i.e. maps of the square which take the edges into $X_1$ and the
vertices into $X_0$, and then take homotopy classes of such maps
relative to the vertices of $I^2$  to form a set $\rho_2 X_*$. Of
course this set will not inherit a group structure but the surprise
is that it does inherit the structure of double groupoid with
connections -- the proof is not entirely trivial, and is given in
\cite{bh1978} and the expository article \cite{brown:hha}. In the
case $X_0$ is a singleton, the equivalence of such double groupoids
to crossed modules takes $\rho X_*$ to the usual second relative
homotopy crossed module.

Thus a search for a {\it higher homotopy groupoid} was  realised in
dimension 2. Connes suggests in \cite{connes} that it has been
fashionable for  mathematicians to disparage groupoids, and it might
be that a lack of attention  to this notion  was one reason why such
a construction had not  been found earlier than 40 years after
Hurewicz's papers.

Finding a good homotopy double groupoid led rather quickly, in view
of the previous experience, to a substantial account of a
2-dimensional HHvKT \cite{bh1978}. This recovers Theorem W, and also
leads to new calculations in 2-dimensional homotopy theory, and in
fact to some new calculations of 2-types. For a recent summary of
some results and some new ones, see the paper in the J. Symbolic
Computation \cite{brownwens:comp} -- publication in this journal
illustrates that we are interested in using general methods in order
to obtain specific calculations, and ones to which there seems no
other route.

Once the 2-dimensional case had been completed in 1975, it was easy
to conjecture the form of general results for dimensions $>2$. These
were proved by 1979 and announcements were made in \cite{bh:CR2}
with full details in \cite{bh:algcub,bh:colimits}. However, these
results needed a number of new ideas, even just to construct the
higher dimensional compositions, and the proof of the HHvKT was
quite hard and intricate. Further, for applications, such as to
explain how the general $\Pi$ behaved on homotopies,  we also needed
a theory of tensor products, found in \cite{bh:tens}, so that the
resulting theory is quite complex. It is also remarkable that ideas
of Whitehead in \cite{jhcw:CHII} played a key role in these results.

\section{Main results}
\label{sec:main} Major features of the work over the years with
Philip Higgins and others can be summarised in the following diagram
of categories and functors:

\begin{diag}
\begin{equation*}
\def\labelstyle{\textstyle}\qquad \vcenter{ \xymatrix@R=3pc{&& \txt{\rm filtered
spaces} \ar @/_2pc/[dll]_(0.6){C_*= \nabla\circ \Pi\;}
  \ar @<-0.5ex>[dl] _-\Pi \ar @<0.5ex>[dr] ^ -\rho & \txt{\rm  filtered \\\rm  cubical sets} \ar [l]_-{\mid \;\mid}\\
\txt{\rm operator \\ \rm chain \\ \rm complexes } \ar @<-0.5ex>
[r]_\Theta  &
\ar @<-0.5ex>[l]_\nabla \txt{\rm crossed \\
\rm complexes} \ar @<-0.5ex>[ur] _-{\mathcal{B}} \ar@<0.5ex>
[rr]^-\lambda  && \ar [u]_{U_*}\txt{\rm cubical \\\rm
$\omega$-groupoids
\\\rm with connections } \ar @<0.5ex>[ll] ^-\gamma   }} %\tag{Main Diagram}
\end{equation*}
\end{diag}in which
\begin{enumerate}
\item[4.1.1] the categories $\FTop$ of filtered spaces, \ogpd\ of
cubical $\omega$-groupoids with connections, and $\crs$ of crossed
complexes are monoidal closed, and have a notion of homotopy using
$\otimes$ and a unit interval object; \item[4.1.2] $\rho, \;\Pi$ are
homotopical functors (that is they are defined in terms of homotopy
classes of certain maps), and preserve homotopies; \item[4.1.3]
  $\lambda, \; \gamma$ are inverse adjoint equivalences
  of monoidal closed categories;
\item[4.1.4]  there is a natural equivalence $\gamma \rho \simeq
\Pi$, so that either $\rho$ or $\Pi$ can be used as appropriate;
\item[4.1.5] $\rho, \;\Pi$ preserve certain colimits
   and certain tensor products;
\item[4.1.6] the category of chain complexes with (a groupoid) of operators
is monoidal closed, $\nabla$ preserves the monoid structures, and is
left adjoint to $\Theta$;
\item[4.1.7] by definition, the {\it cubical filtered classifying
space} is $\mathcal{B^\Box}=|\,|\circ U_*$ where $U_*$ is the
forgetful functor to filtered cubical sets\footnote{Cubical sets are
defined, analogously to simplicial sets, as functors $K:
\square^{op} \to \Set$ where $\square$ is the `box' category with
objects $I^n$ and morphisms the compositions of inclusions of faces
and of the various projections $I^n \to I^r$ for $n>r$. The {\it
geometric realisation} $|K|$ of such a cubical set is obtained by
quotienting the disjoint union of the sets $K(I^n)\times I^n$ by the
relations defined by the morphisms of $\square$. For more details,
see \cite{Jardine-cat-cub}, and for variations on the category
$\square$ to include for example connections, see
\cite{grandismauri}. See also section \ref{sec:omegagroupoid}.}
using the filtration of an $\omega$-groupoid by skeleta, and $|\,|$
is geometric realisation of a cubical set;
 \item[4.1.8] there is a natural equivalence $\Pi
 \circ\mathcal{B^\Box}\simeq 1$;
\item[4.1.9] \label{blank:homclass} if $C$ is a crossed complex and its cubical
classifying space is defined as $B^\Box
C=(\mathcal{B^\Box}C)_\infty$, then for a $CW$-complex $X$, and
using homotopy as in 4.1.1 for crossed complexes, there is a natural
bijection of sets of homotopy classes
$$[X,B^\Box C] \cong [\Pi X_*,C]. $$
Recent applications of the simplicial version of the classifying space  are
in \cite{brown:homclass,Porter-Tur,martinsporter}.
\end{enumerate}

Here a {\it filtered space} consists of a (compactly generated)
space $X_\infty$ and an increasing sequence of subspaces
$$X_*: X_0 \subset X_1 \subset X_2 \subset \cdots \subset
X_\infty. $$ With the obvious morphisms, this gives the category
$\FTop$. The tensor product in this category is the usual $$(X_*
\otimes Y_*)_n = \bigcup _{p+q=n} X_p \times Y_q.$$ The closed
structure is easy to construct from the law $$ \FTop(X_* \otimes
Y_*, Z_*) \cong \FTop(X_*, \FTOP(Y_*,Z_*)).$$ An advantage of this
monoidal closed structure is that it allows an enrichment of the
category $\FTop$ over either crossed complexes or \ogpd\ using
$\Pi$ or $\rho$ applied to $\FTOP(Y_*,Z_*)$.

The structure of {\it crossed complex} is suggested by the canonical
example, the {\it fundamental crossed complex} $\Pi X_*$ of the
filtered space $X_*$. So it is given by a diagram \begin{diag}
 $$\diagram \cdots \rto & C_n \dto<-.05ex>^{t}\rto^-{\delta_n}
& C_{n-1} \rto \dto<-1.2ex>^{t} & \cdots \rto &
 C_2\rto^-{\delta_2} \dto<-.05ex>^{t} & C_1
\dto<0.0ex>^(0.45){t} \dto<-1ex>_(0.45){s} \\ &
C_0&C_0\rule{0.5em}{0ex}  & & \rule{0.5em}{0ex} C_0 &
\rule{0em}{0ex} C_0  \enddiagram$$
\end{diag}
in which in this example  $C_1$ is the fundamental groupoid $
\pi_1(X_1,X_0)$ of $X_1$ on the `set of base points' $C_0= X_0$,
while for $n \ge 2$,   $C_n$ is the family of relative homotopy
groups $\{ C_n(x) \} =\{\pi_n(X_n,X_{n-1},x) \, | \,x \in X_0 \}$.
The boundary maps are those standard in homotopy theory. There is
for $n \ge 2$ an action of the groupoid $C_1$ on $C_n$ (and of
$C_1$ on the groups $C_1(x), \, x \in X_0$ by conjugation), the
boundary morphisms are operator morphisms,
$\delta_{n-1}\delta_n=0, \, n \ge 3$, and the additional axioms
are satisfied that
\begin{blank} $b\io   c b = c^{\delta_2 b}, \, b,c \in C_2$,
so that $\delta_2:  C_2 \to C_1$ is a crossed module (of groupoids);
\end{blank}
\begin{blank}\label{blank:trivact}  if $c \in C_2$ then $\delta_2 c$ acts trivially on $C_n$
for $n \ge 3$;
\end{blank}
\begin{blank} each group $C_n(x)$ is
Abelian for $n \ge 3$, and so the family $C_n$ is a $C_1$-module.
\end{blank}
Clearly we  obtain a category $\crs$ of crossed complexes; this
category is not so familiar and so we give arguments for using it
in the next section.

As algebraic examples of crossed complexes we have: $C=\bC(G,n)$
where $G$ is a group, commutative if $n \ge 2$, and $C$ is $G$ in
dimension $n$ and trivial elsewhere; $C=\bC(G,1:M,n)$, where $G$ is
a group, $M$ is a $G$-module, $n\ge 2$, and $C$ is $G$ in dimension
1, $M$ in dimension $n$, trivial elsewhere, and zero boundary if
$n=2$; $C$ is a crossed module (of groups) in dimensions 1 and 2 and
trivial elsewhere.

A crossed complex $C$ has a fundamental groupoid $\pi_1 C =
C_1/\Im \delta_2$, and also for $n  \ge 2$ a family $\{H_n(C,p)|p
\in C_0\}$ of homology groups.

\section{Why crossed complexes?}
\label{sec:whycrossed}

\hspace{1em} $\bullet$ They generalise groupoids and crossed
modules to all dimensions. Note that the natural context for
second relative homotopy groups is crossed modules of groupoids,
rather than groups.

$\bullet$ They are good for modelling $CW$-complexes.

$\bullet$ Free crossed resolutions enable calculations with small
$CW$-complexes and  $CW$-maps, see section \ref{S:freecrossed}.

$\bullet$ Crossed complexes  give a kind of `linear model' of
homotopy types which includes all 2-types. Thus although they are
not the most general model by any means (they do not contain
quadratic information such as Whitehead products), this simplicity
makes them easier to handle and to relate to classical tools. The
new methods and results obtained for crossed complexes can be used
as a model for more complicated situations. This is how a general
$n$-adic Hurewicz Theorem was found \cite{brlo:hur}.

$\bullet$ They are convenient for {\it calculation}, and the
functor $\Pi$ is classical, involving relative homotopy groups. We
explain some results in this form later.

$\bullet$ They are close to chain complexes with a group(oid) of
operators, and related to some classical homological algebra (e.g.
{\it chains of syzygies}). In fact if $SX$ is the simplicial
singular complex of a space, with its skeletal filtration, then
the crossed complex $\Pi(SX)$ can be considered as a slightly non
commutative version of the singular chains of a space.

$\bullet$ The monoidal structure is suggestive of further
developments (e.g. {\it crossed differential algebras}) see
\cite{bauesT1,babr}. It is used in \cite{br:gilbert} to give an
algebraic model of homotopy $3$-types, and to discuss
automorphisms of crossed modules.

$\bullet$ Crossed complexes have a good homotopy theory, with a {\it
cylinder object, and  homotopy colimits}, \cite{br-go}. The homotopy
classification result 4.1.9 generalises a classical theorem of
Eilenberg-Mac~Lane. Applications of (the simplicial version) are
given in for example \cite{martins, martinsporter,Porter-Tur}.

$\bullet$ They have an interesting relation with the Moore complex
of simplicial groups and of simplicial groupoids (see section
\ref{S:simplicial}).

\section{Why cubical $\omega$-groupoids with  connections?}
\label{sec:whycubomega}

The definition of these objects is more difficult to give, but will
be indicated in section \ref{sec:omegagroupoid}. Here we explain why
these structures are a kind of engine giving the power behind the
theory.

$\bullet$ The functor $\rho$ gives a form of {\it higher homotopy
groupoid}, thus confirming the visions of the early topologists.

$\bullet$ They are equivalent to crossed complexes.

$\bullet$ They have a clear {\it monoidal closed structure}, and a
notion of homotopy, from which one can deduce those on crossed
complexes, using the equivalence of categories.

$\bullet$ It is easy to relate the functor $\rho$ to tensor
products, but quite difficult to do this directly for $\Pi$.

$\bullet$ Cubical   methods,  unlike  globular or simplicial
methods, allow for a simple  {\it algebraic inverse to
subdivision}, which is crucial for our local-to-global theorems.

$\bullet$ The additional structure of `connections',  and the
equivalence with crossed complexes,  allows for the sophisticated
notion of {\it commutative cube}, and the proof that {\it multiple
compositions of commutative cubes are commutative}. The last fact is
a key component of the proof of the HHvKT.

$\bullet$ They yield a construction  of a {\it (cubical)
classifying space} $B^\Box C=(\mathcal{B^\Box}C)_\infty$ of a
crossed complex $C$, which generalises (cubical) versions of
Eilenberg-Mac~Lane spaces, including the local coefficient case.
This has convenient relation to homotopies.

$\bullet$ There is a current {\it resurgence of the use of cubes} in
for example combinatorics, algebraic topology, and concurrency.
There is a Dold-Kan type theorem for cubical Abelian  groups with
connections \cite{bh:cubabgp3}.

\section{The equivalence of categories}
\label{sec:equivcat}

Let $\crs$, \ogpd\ denote respectively the categories of crossed
complexes and $\omega$-groupoids: we use the latter term as an
abbreviation of `cubical $\omega$-groupoids with connections'. A
major part of the work consists in defining these categories and
proving their equivalence, which thus gives an example of two
algebraically defined categories whose equivalence is non trivial.
It is even more subtle than that because the functors $\gamma: \crs
\to \ogpdm$, $\lambda : \ogpdm \to \crs$ are not hard to define, and
it is easy to prove $\gamma \lambda \simeq 1$. The hard part is to
prove $ \lambda \gamma \simeq 1$, which  shows that an
$\omega$-groupoid $G$ may be reconstructed from the crossed complex
$\gamma(G)$ it contains. The proof involves using the connections to
construct a `folding map' $\Phi: G_n \to G_n$ , with image
$\gamma(G)_n$, and establishing its major properties, including the
relations with the compositions. This gives an algebraic form of
some old intuitions of several ways of defining relative homotopy
groups, for example using cubes or cells.

On the way we establish properties of {\it thin elements}, as those
which fold down to 1, and show that $G$ satisfies a strong Kan
extension  condition, namely that every box has a unique thin
filler. This result plays a key role in the proof of the HHvKT for
$\rho$, since it is used to show an independence of choice. That
part of the proof goes by showing that the two choices can be seen,
since we start with a homotopy, as given by the two ends $\partial
^\pm_{n+1} x$ of an $(n+1)$-cube $x$. It is then shown by induction,
using the method of construction and the above result, that $x$ is
degenerate in direction $n+1$. Hence the two ends in that direction
coincide.

Properties of the folding map are used also in showing that $\Pi X_*
$ is actually included in $\rho X_* $; in relating two types of
thinness for elements of $\rho X_* $; and in proving a {\it homotopy
addition lemma} in $\rho X_* $.

Any \ogpd\ $G$ has an underlying cubical set $UG$. If $C$ is a
crossed complex, then the cubical set  $U(\lambda C)$ is called
the {\it cubical nerve } $N^\Box C$ of $C$. It is a conclusion of
the theory that we can also obtain $N^\Box C$ as
$$(N^\Box C)_n= \crs(\Pi I^n_*, C)$$ where $I^n_*$ is the usual
geometric cube with its standard skeletal filtration. The (cubical)
geometric realisation $|N^\Box C|$ is also called the {\it cubical
classifying space} $B^\Box C$ of the crossed complex $C$. The
filtration $C^*$ of $C$ by skeleta gives a filtration $B^\Box C^*$
of $B^\Box C$ and there is (as in 4.1.6) a natural isomorphism
$\Pi(B^\Box C^*) \cong C$. Thus the properties of a crossed complex
are those that are universally satisfied by $\Pi X_* $. These proofs
use the equivalence of the homotopy categories of Kan\footnote{The
notion of Kan cubical set $K$ is also called a cofibrant cubical
set. It is an extension condition that any partial $r$-box in $K$ is
the partial boundary of an element of $K_r$. See for example
\cite{Jardine-cat-cub}, but the idea goes back to the first paper by
D. Kan in 1958.  } cubical sets and of $CW$-complexes. We originally
took this from the Warwick Masters thesis of S. Hintze, but it is
now available with different proofs from Antolini
\cite{antolini:thesis} and Jardine \cite{Jardine-cat-cub}.

As said above, by taking particular values for $C$, the
classifying space $B^\Box C$ gives cubical versions of
Eilenberg-Mac~Lane spaces $K(G,n)$, including the case $n=1$ and
$G$ non commutative. If $C$ is essentially a crossed module, then
$B^\Box C$ is called the {\it cubical classifying space} of the
crossed module, and in fact realises the $k$-invariant of the
crossed module.

Another useful result is that if $K$ is a cubical set, then
$\rho(|K|_*)$ may be identified with $\rho(K)$, the {\it free \ogpd\
on the cubical set $K$}, where here $|K|_*$ is the usual filtration
by skeleta. On the other hand, our proof that $\Pi(|K|_*)$ is the
free crossed complex on the non-degenerate cubes of $K$ uses the
generalised HHvKT of the next section.

It is also possible to give simplicial and globular
 versions of
some of the above results, because the category of crossed
complexes is equivalent also to those of simplicial $T$-complexes
\cite{ashley} and of globular $\infty$-groupoids \cite{bh1981a}.
In fact the published paper on the classifying space of a crossed
complex \cite{bh:class91} is given in simplicial terms, in order
to link more easily with well known theories.

\section{First main aim of the work: Higher Homotopy van Kampen Theorems}
\label{sec:HHvKT} These theorems give {\it non commutative  tools
for higher dimensional local-to-global problems} yielding  a variety
of new, often  non commutative, calculations, which {\it prove}
(i.e. test) the theory. We now explain these theorems in a way which
strengthens the relation with descent, since that was a theme of the
conference at which the talk was given on which this survey is
based.

We suppose given an open cover $\cU = \{ U^\lambda \}_{\lambda \in
\Lambda}$ of $X$. This cover defines a map
$$q: E=\bigsqcup_{\lambda \in \Lambda}U^\lambda \to X$$ and
so we can form an augmented  simplicial space
$$\Ce(q): \xymatrix{\cdots  \ar@<1ex>[r] \ar [r] \ar@<-1ex>[r]E \times_{X} E \times_{X}
E & E \times_{X} E \ar@<0.5ex>[r]  \ar@<-0.5ex>[r]& E \ar [r]^q &
X } $$ where the higher dimensional terms involve disjoint unions
of multiple intersections $U^\nu$ of the $U^\lambda$.

We now suppose given a filtered space  $X_*$, a cover $\cU$ as
above of $X=X_\infty$, and so an augmented simplicial filtered
space $\Ce(q_*)$ involving multiple intersections  $U^\nu_*$ of
the induced filtered spaces.

We still need a connectivity condition.

\begin{Def} A filtered space $X_*$ is {\it connected} if and only if
 the induced maps $\pi_0 X_0 \rightarrow \pi_0 X_n$ are
surjective and  $\pi_n (X_r , X_n , \nu ) = 0$  for all $ n
> 0, r > n $ and $\nu \in X_0$.  \end{Def}

\begin{thm}[MAIN RESULT (HHvKT)]\label{thm:HHvKT} If $U^\nu_*$ is connected for all
finite  intersections $U^\nu$ of the elements of the open cover, then\\
{\em (C) (connectivity)} $X_*$ is connected, and \\{\em (I)
(isomorphism)}  the following diagram as part of $\rho(\Ce(q_*))$
\begin{equation}\label{rhodiag}  \xymatrix{ \rho(E_* \times_{X_*}
E_*)\ar@<0.5ex>[r] \ar@<-0.5ex>[r] & \rho E_*\ar [r]^{\rho(q_*)}&
\rho X_* .}\tag*{ {(c$\rho$)}}\end{equation} is a coequaliser
diagram. Hence the following diagram of crossed complexes
 {\begin{equation} \xymatrix{ \Pi(E_* \times_{X_*}
E_*)\ar@<0.5ex>[r] \ar@<-0.5ex>[r] & \Pi E_*\ar [r]^{\Pi(q_*)}& \Pi
X_* .}\tag*{ {(c$\Pi$)}}\end{equation}}is also a   coequaliser
diagram.
 \end{thm}
So we get calculations of the fundamental crossed complex $\Pi X_*
$.

It should be emphasised that to get to and apply this theorem takes
just the two papers \cite{bh:algcub,bh:colimits} totalling   58
pages. With this we deduce in the first instance:
\begin{itemize}
  \item the usual vKT for the fundamental groupoid on a set of
  base points;
  \item  the Brouwer degree theorem  ($\pi _n S^n=\Z$);
  \item   the relative  Hurewicz theorem;
  \item    Whitehead's theorem that  $\pi _n(X \cup
\{e^2_{\lambda} \},X)$  is a free crossed module;
    \item  an excision result, more general than the previous two,  on $\pi_n(A \cup B, A,x)$
    as an induced module (crossed module if $n=2$) when  $(A,A\cap
    B)$ is $(n-1)$-connected.
  \end{itemize}
The assumptions required of the reader are quite small, just some
familiarity with $CW$-complexes. This contrasts with some
expositions of basic homotopy theory, where the proof of say the
relative Hurewicz theorem requires knowledge of singular homology
theory. Of course it is surprising to get this last theorem without
homology, but this is because it is seen as a statement on the
morphism of relative homotopy groups $$\pi_n(X,A,x) \to \pi_n(X \cup
CA,CA,x)\cong \pi_n(X \cup CA,x)$$ and is obtained, like our proof
of Theorem W, as a special case of an excision result. The reason
for this success is that {\it we use algebraic structures which
model the underlying processes of the geometry}  more closely than
those in common use. These algebraic structures and their relations
are quite intricate, as befits the complications of homotopy theory,
so the theory is tight knit.

Note also that these results cope well with the action of the
fundamental group on higher homotopy groups.

The calculational use of the HHvKT for $\Pi X_* $ is enhanced by the
relation of $\Pi$  with tensor products (see section
\ref{S:freecrossed} for more details).

\section{The fundamental
cubical  \om-groupoid $\rho X_* $  of a filtered space $X_*$}
\label{sec:omegagroupoid}
 Here are the basic elements of the construction.

$I^n_*$: the $n$-cube with its skeletal filtration.

Set  $R_n X_* = \FTop(I^n_*, X_*)$.  This is a  {\it cubical set
with compositions, connections, and inversions}.

For $i=1, \ldots,n$ there are standard:

  {{\em face maps}} $\pt^\pm_i: R_nX_*
\to R_{n-1}X_*$;

  {{\em degeneracy maps}} $\eps_i :R_{n-1}X_* \to
R_{n}X_*$

  {{\em connections}} $ \Gamma_i^\pm :R_{n-1 }X _*\to R_nX_*$

    {{\em compositions}} $a\circ_ib$ defined for $a,b\in R_nX_*$ such
that $\pt^+_ia=\pt^-_ib$

  { {\em inversions}} $-_i: R_n \to R_n$.

The connections are induced  by $\gamma_i^\A : I^n \to
 I^{n-1}$    defined using the monoid structures $\max,\min :I^2 \to
 I$. They are essential for many reasons, e.g. to discuss the notion of
  {\it commutative cube}.

These operations have certain algebraic properties which are easily
derived from the geometry and  which we do not itemise here -- see
for example \cite{aa-b-s}. These were listed first in the Bangor
thesis of Al-Agl \cite{alagl:thesis}. (In the paper \cite{bh:algcub}
the only basic connections needed are the $\Gamma^+_i$, from which
the $\Gamma^-_i$ are derived using the inverses of the groupoid
structures.)

Now it is  natural and convenient  to define $f \equiv g$ for $f,g:
I^n_* \to X_*$ to mean $f$ is homotopic to $g$ through filtered maps
an relative to the vertices of $I^n$. This gives a quotient map
{$$p: R_n X_* \to \rho_n X_* =(R_n X_* /\equiv). $$}

The following results are proved in \cite{bh:colimits}.

\begin{blank} \label{comp-inher}
The compositions on $RX_*$ are inherited by $\rho X_*$ to give $\rho
X_* $ the structure of cubical multiple groupoid with connections.
\end{blank}

\begin{blank}The map $p: R X_* \to \rho X_* $ is a Kan fibration of
cubical sets. \label{fibr}
\end{blank}

The proofs of both  results use methods of collapsing which are
indicated in the next section. The second result is almost
unbelievable. Its proof has to give a systematic method of deforming
a cube with the right faces `up to homotopy' into a cube with
exactly the right faces, using the given homotopies. In both cases,
the assumption that the relation $\equiv$ uses homotopies relative
to the  vertices is essential to start the induction. (In fact the
paper \cite{bh:colimits} does not use homotopy relative to the
vertices, but imposes an extra condition $J_0$, that each loop in
$X_0$ is contractible $X_1$, which again starts the induction. This
condition is awkward in applications, for example to function
spaces. A full exposition of the whole story is in preparation,
\cite{bhs:nonabalgtop}.)

An essential ingredient in the proof of the HHvKT is the notion of
{\it multiple composition}. We have discussed this already in
dimension 2, with a suggestive picture in the diagram
(\ref{eq:multcomp}).  In dimension $n$, the aim is to give algebraic
expression to the idea of a cube $I^n$ being subdivided by
hyperplanes parallel to the faces into many small cubes, a
subdivision with a long history in mathematics.

Let $(m) = (m_1, \ldots , m_n)$ be an $n$-tuple of positive integers
and
$$\phi_{(m)} : I^n \rightarrow [0, m_1] \times \cdots \times [0,
m_n]$$ be the map $(x_1 , \ldots , x_n) \mapsto (m_1 x_1, \ldots ,
m_n x_n).$ Then a {\it subdivision of type $(m)$}  of a map $\alpha
: I^n \rightarrow X$ is a factorisation $\alpha = \alpha' \circ
\phi_{(m)}$; its {\it parts} are the cubes $\alpha_{(r)}$ where $(r)
= (r_1, \ldots , r_n)$ is an $n$-tuple of integers with $1 \leqslant
r_i \leqslant m_i$, $i = 1, \ldots , n,$ and where $\alpha_{(r)} :
I^n \rightarrow X$ is given by
$$(x_1, \ldots , x_n) \mapsto
\alpha'(x_1 + r_1 - 1, \ldots , x_n + r_n - 1).$$

We then say that $\alpha$ is the {\it composite}  of the cubes
$\alpha_{(r)}$ and write $\alpha = [\alpha_{(r)}]$. The {\it domain}
of $\alpha_{(r)}$ is then the set $\{(x_1,\ldots,x_n) \in I^n :
r_i-1 \leqslant x_i \leqslant r_i, 1 \leqslant i \leqslant n\}$.
This ability to express `algebraic inverse to subdivision'  is one
benefit of using cubical methods.

Similarly, in a cubical set with compositions satisfying the
interchange law we can define the multiple composition
$[\alpha_{(r)}] $ of a multiple array $(\alpha_{(r)})$ provided the
obviously necessary multiple incidence relations  of the individual
$\alpha_{(r)}$ to their neighbours are satisfied.

Here is an application which is essential in many proofs, and which
seems hard to prove without the techniques involved in \ref{fibr}.

\begin{thm}[Lifting multiple compositions]
Let  $[\alpha_{(r)}]$ be a multiple composition in $\rho_n X_* $.
Then representatives  $a_{(r)}$ of the $\alpha_{(r)}$   may be
chosen so that the multiple composition $[a_{(r)}]$ is well defined
in $R_n X_* $.
\end{thm}
\noindent {\bf  Proof:} The multiple composition $[\alpha_{(r)}]$
determines a cubical map
$$A:K \to\rho X_* $$
where the cubical set $K$ corresponds to  a representation of the
multiple composition by a subdivision of the geometric cube, so
that top cells $c_{(r)}$ of $K$ are mapped by $A$ to
$\alpha_{(r)}$.

Consider the diagram, in which $*$ is a corner vertex of $K$,
$$\xymatrix@=4pc{\ast \ar [r] \ar [d] & R X_*  \ar [d] ^p \\
K \ar [r] _-A \ar @{-->}[ur]^{A'} & \rho X_* }.$$ Then $K$ collapses
to $*$, written $K\searrow *$. (As an example, see how the
subdivision in the diagram (\ref{eq:multcomp}) may be collapsed row
by row to a point.)  By the fibration result, $A$ lifts to $A'$,
which represents $[a_{(r)}]$, as required. \rule{0ex}{1mm}\hfill
$\Box$

So we have to explain collapsing.

\section{Collapsing}
\label{sec:collapsing} We use a basic notion of collapsing and
expanding due to J.H.C. Whitehead, \cite{wjhc:sht}.

 Let $C
\subset B$ be subcomplexes of $I^n $. We say $C$ is an  {\it
elementary collapse} of $B$,  $B\ecoll C,$ if for some $s \geq 1$
there is an $s$-cell $a$ of $B$ and $(s - 1)$-face $b$ of $a$, the
{\em free face},  such that $$B = C \cup a, \qquad C \cap a =
\dot{a} \setminus b$$ (where $\dot{a}$ denotes the union of the
proper faces of $a$).

We say  $B_1$  {\it collapses} to $B_r$, written $B_1 \searrow B_r$,
if there is a sequence
$$B_1 \ecoll B_2 \ecoll \cdots \ecoll B_r$$ of elementary collapses.

If $C$ is a subcomplex of $B$ then $$B \times I \searrow  (B
\times \{ 0 \} \cup C \times I)$$ (this is proved by induction on
dimension of $B \setminus C$).

Further, $I^n$ collapses to any one of its vertices (this may be
proved by induction on $n$ using the first example). These
collapsing techniques allows the construction of the extensions of
filtered maps and filtered homotopies that are crucial for proving
\ref{comp-inher}, that $\rho X_* $ does obtain the structure of
 multiple groupoid.

However,  more subtle collapsing techniques using partial boxes
are required to prove the fibration theorem \ref{fibr},  as partly
explained in the next section.

%\end{example}

\vspace{1cm}

\section{Partial boxes}
\label{sec:partialbox} Let $C$ be an $r$-cell in the $n$-cube $I^n.$
Two $(r - 1)$-faces of $C$ are called  {\it opposite} if they do not
meet.

A  {partial box} in $C$ is a subcomplex $B$ of $C$ generated by
one $(r - 1)$-face $b$ of $C$ (called a {\it base} of $B$) and a
number, possibly zero, of other $(r - 1)$-faces of $C$ none of
which is opposite to $b.$

The partial box is a  {box} if its $(r - 1)$-cells consist of all
but one of the $(r - 1)$-faces of $C.$

The proof of the fibration theorem uses a filter homotopy
extension property and the following:

\begin{prop}[Key Proposition] Let $B, B'$ be partial boxes in an
$r$-cell $C$ of $I^n$ such that $B' \subset B.$ Then there is a
chain
$$B = B_s \searrow B_{s-1} \searrow \cdots \searrow B_1 = B'$$
such that \begin{enumerate}[ \rm (i)]
\item each $B_i$ is a partial box in $C$;
\item $B_{i+1} = B_i \cup a_i$ where $a_i$ is an $(r -
 1)$-cell of $C$ not in $B_i$;
 \item $a_i \cap B_i$ is a
 partial box in $a_i.$

\end{enumerate}
\end{prop}
The proof is quite neat, and follows the pictures. Induction up such
a chain of partial boxes is one of the steps in the proof of the
fibration theorem \ref{fibr}. The proposition implies that an
inclusion of partial boxes is what is known as an anodyne extension,
\cite{Jardine-cat-cub}.

Here is an example of a sequence of collapsings of a partial box
$B$, which illustrate some choices in forming a collapse  $B
\searrow {\mathbf 0}$ through two other partial boxes $B_1, B_2$.

\vspace{4mm}

\setlength{\unitlength}{0.00033333in}
\begingroup\makeatletter\ifx\SetFigFont\undefined
% extract first six characters in \fmtname
\def\x#1#2#3#4#5#6#7\relax{\def\x{#1#2#3#4#5#6}}%
\expandafter\x\fmtname xxxxxx\relax \def\y{splain}%
\ifx\x\y   % LaTeX or SliTeX?
\gdef\SetFigFont#1#2#3{%
  \ifnum #1<17\tiny\else \ifnum #1<20\small\else
  \ifnum #1<24\normalsize\else \ifnum #1<29\large\else
  \ifnum #1<34\Large\else \ifnum #1<41\LARGE\else
     \huge\fi\fi\fi\fi\fi\fi
  \csname #3\endcsname}%
\else \gdef\SetFigFont#1#2#3{\begingroup
  \count@#1\relax \ifnum 25<\count@\count@25\fi
  \def\x{\endgroup\@setsize\SetFigFont{#2pt}}%
  \expandafter\x
    \csname \romannumeral\the\count@ pt\expandafter\endcsname
    \csname @\romannumeral\the\count@ pt\endcsname
  \csname #3\endcsname}%
\fi \fi\endgroup {\renewcommand{\dashlinestretch}{30}
\begin{picture}(6777,1964)(1000,-10)
\drawline(400,12)(1580,12)(2000,715)(820,715)(400,12)
\drawline(400,12)(400,1190) \drawline(820,715)(820,1893)
\drawline(820,1893)(400,1190) \drawline(1580,12)(1580,1190)
\drawline(2000,715)(2000,1893) \drawline(1580,1190)(2000,1893)

\put(1000,-500){\makebox(0,0)[lb]{{\SetFigFont{10}{14.4}{rm}$B$}}}
\put(9100,-500){\makebox(0,0)[lb]{{\SetFigFont{10}{14.4}{rm}$B_1$}}}
\put(14000,-500){\makebox(0,0)[lb]{{\SetFigFont{10}{14.4}{rm}$B_2$}}}

\drawline(2300,1174)(2700,634)
\drawline(2618,696)(2700,634)(2657,727)
\put(2550,980){\makebox(0,0)[lb]{{\SetFigFont{10}{14.4}{rm}$e$}}}

\drawline(3200,12)(4380,12)(4800,715)(3620,715)(3200,12)
\drawline(3620,715)(3620,1893) \drawline(3200,12)(3200,1190)
\drawline(3620,1893)(3200,1190) \drawline(4380,12)(4380,1190)
\drawline(4800,715)(4800,1893)

\drawline(5100,1174)(5500,634)
\drawline(5418,696)(5500,634)(5457,727)
\put(5350,980){\makebox(0,0)[lb]{{\SetFigFont{10}{14.4}{rm}$e$}}}

\drawline(6000,12)(7180,12)(7600,715)(6420,715)(6000,12)
\drawline(6420,715)(6420,1893) \drawline(6000,12)(6000,1190)
\drawline(6420,1893)(6000,1190) \drawline(7600,715)(7600,1893)

\drawline(7900,1174)(8300,634)
\drawline(8198,696)(8300,634)(8257,727)
\put(8150,980){\makebox(0,0)[lb]{{\SetFigFont{10}{14.4}{rm}$e$}}}

\drawline(8800,12)(9980,12)(10400,715)(9220,715)(8800,12)
\drawline(9220,715)(9220,1893) \drawline(8800,12)(8800,1190)
\drawline(9220,1893)(8800,1190)

\drawline(10700,1174)(11100,634)
\drawline(11018,696)(11100,634)(11057,727)
\put(10950,980){\makebox(0,0)[lb]{{\SetFigFont{10}{14.4}{rm}$e$}}}

\put(11800,437){\makebox(0,0)[lb]{{\SetFigFont{10}{14.4}{rm}$\cdots$}}}

\drawline(12700,1174)(13100,634)
\drawline(13018,696)(13100,634)(13057,727)
\put(12950,980){\makebox(0,0)[lb]{{\SetFigFont{10}{14.4}{rm}$e$}}}

\drawline(14300,715)(13880,12) \drawline(14300,715)(14300,1893)
\drawline(13880,12)(13880,1190) \drawline(14300,1893)(13880,1190)

\drawline(14800,1174)(15200,634)
\drawline(15118,696)(15200,634)(15157,727)
\put(15050,980){\makebox(0,0)[lb]{{\SetFigFont{10}{14.4}{rm}$e$}}}

\put(15900,437){\makebox(0,0)[lb]{{\SetFigFont{10}{14.4}{rm}$\cdots$}}}

\drawline(16800,1174)(17200,634)
\drawline(17118,696)(17200,634)(17157,727)
\put(17050,980){\makebox(0,0)[lb]{{\SetFigFont{10}{14.4}{rm}$e$}}}

\put(18000,715){\makebox(0,0)[lb]{{\SetFigFont{14}{14.4}{rm}$\cdot$}}}
\end{picture}
} \vspace{0.75cm}

The proof of the fibration theorem  gives a program for carrying
out the deformations needed to do the lifting. In some sense, it
implies  computing a multiple composition can be done using
collapsing as the guide.

Methods of collapsing generalise methods of trees in dimension 1.

\section{Thin elements}
\label{sec:thinelements}

Another key concept is that of  {\it thin element} $\alpha \in
\rho_n X_* $ for $n\ge 2$. The proofs here use strongly  results of
\cite{bh:algcub}.

We say $\alpha$ is  {\it geometrically thin} if it has a {\it
deficient} representative, i.e. an  $a: I^n_* \to X_*$ such that
$a(I^n) \subset X_{n-1}$.

We say $\alpha$ is  {\it algebraically thin} if it is a multiple
composition of degenerate elements or those coming from repeated
(including 0) negatives of connections. Clearly any multiple
composition of algebraically thin elements is thin.

\begin{thm}  {\rm (i)} Algebraically thin is equivalent to geometrically thin.

{\rm (ii)} In a cubical \om-groupoid with connections, any box has a
unique thin filler.
\end{thm}
\begin{proof} The proof of the forward implication in (i) uses lifting of multiple
compositions, in a stronger form than stated above.

The proofs of (ii) and the backward implication in (i) use the full
force of the algebraic relation between \om-groupoids and crossed
complexes.
\end{proof}

These results allow one to replace arguments with commutative
cubes by arguments with thin elements.

\section{Sketch proof of the HHvKT}
\label{sec:sketchproofHHvKt} The proof goes by verifying the
required universal property. Let $\mathcal{U}$ be an open cover of
$X$ as in Theorem \ref{thm:HHvKT}.

We go back to the following diagram whose top row is  part of
$\rho(\Ce(q_*))$ {
\begin{equation} \xymatrix@R=3pc{ \rho(E_* \times_{X_*}
E_*)\ar@<0.5ex>[rr]^-{\partial_0} \ar@<-0.5ex>[rr]_-{\partial_1} &&
\ar [drr] _f \rho(E_*)\ar [rr]^{\rho(q_*)}&& \rho X_*  \ar
@{-->}[d]^{f'}\\&&&&G }\tag*{ {(c$\rho$)}}\end{equation}} To prove
this top row is a coequaliser diagram, we suppose given a morphism
$f:\rho(E_*) \to G$ of cubical $\omega$-groupoids with connection
such that $f \circ \partial_0=f \circ \partial_1$, and prove that
there is a unique morphism $f':\rho X_*  \to G$ such that
$f'\circ\rho(q_*)=f$.

To define $f'(\alpha)$ for $\alpha \in \rho X_* $, you subdivide a
representative $a$ of  $\alpha$ to give  $a = [a_{(r)}]$ so that
each $a_{(r)}$ lies in an element $U^{(r)}$ of $\cU$; use the
connectivity conditions and this subdivision to deform $a$ into
$b=[b_{(r)}]$ so that
\begin{align*} b_{(r)} &\in R(U^{(r)}_*)
\\ \intertext{ and so obtain} \beta_{(r)} &\in \rho(U^{(r)}_*).
\intertext{ The elements}f\beta_{(r)} &\in G \\ \intertext{ may be
composed in $G$ (by the conditions on $f$), to give an element}
\theta(\alpha)=[f\beta_{(r)} ]&\in G.\end{align*}   So the proof
of the universal property has to use an  {\it algebraic inverse to
subdivision}. Again an analogy here is with sending an email: the
element you start with is subdivided, deformed so that each part
is correctly labelled, the separate parts are sent, and then
recombined.

The proof that  $\theta(\alpha)$ is independent of the choices made
uses crucially properties of thin elements. The key point is: {\it a
filter homotopy $h: \alpha \equiv \alpha'  $  in $R_n X_* $ gives a
deficient element of $R_{n+1} X_* $.}

The method is to do the subdivision and deformation argument on
such a homotopy, push the little bits in some
$$\rho_{n+1}(U^\lambda_*)$$ (now thin) over to $G$, combine them
and get a thin element $$\tau \in G_{n+1}$$ all of whose faces not
involving the direction $(n+1)$ are thin  { \it because $h$ was
given to be a filter homotopy}. An inductive argument on unique
thin fillers of boxes then shows that {\it $\tau$ is degenerate in
direction $(n+1)$}, so that the two ends in direction $(n+1)$ are
the same.

This ends a rough sketch of the proof of the HHvKT for $\rho$.

Note that the theory of these forms of multiple groupoids is
designed to make this last argument work. We replace a formula for
saying a cube $h$ has commutative boundary by a statement that $h$
is thin. It would be very difficult to replace the above argument,
on the composition  of thin elements, by a higher dimensional
manipulation of formulae such as that given in section 3 for a
commutative $3$-cube.

Further, the proof does not require knowledge of the existence of
all coequalisers, not does it give a recipe for constructing these
in specific examples.

\section{Tensor products and homotopies}
\label{sec:tenshomotopies}

The construction of the  monoidal closed structure on the category
\ogpd\ is based on rather formal properties of cubical sets, and
the fact that for the cubical set $\bI^n$ we have $\bI^m \otimes
\bI^n \cong \bI^{m+n}$. The details are given in \cite{bh:tens}.
The equivalence of categories implies then that the category
$\crs$ is also monoidal closed, with a natural isomorphism
$$\crs(A \otimes B,C) \cong \crs(A, \Crs(B,C)).$$ Here the
elements of $\Crs(B,C)$ are in dimension 0 the morphisms $B \to
C$, in dimension 1  the {\it left homotopies of morphisms}, and in
higher dimensions are forms of higher homotopies. The precise
description of these is obtained of course by tracing out in
detail the equivalence of categories. It should be emphasised that
certain choices are made in constructing this equivalence, and
these choices are reflected in the final formulae that are
obtained.

An important result is that if $X_*,Y_*$ are filtered spaces, then
there is a natural transformation
\begin{alignat*}{2} \label{eta}
\eta \;:\; \rho X_*  &\otimes \rho Y_* &&\to \rho(X_*\otimes Y_*)\\
                    [a] &\otimes [b]   && \mapsto [a \otimes b]
\end{alignat*}
where if $a:I^m_* \to X_*, \, b: I^n_* \to Y_*$ then $a \otimes b
: I^{m+n}_* \to X_* \otimes Y_*$. It not hard to see, in this
cubical setting, that $\eta$ is well defined.  It can also be
shown using previous results that $\eta$ is an isomorphism if
$X_*,Y_*$ are the geometric realisations of cubical sets with the
usual skeletal filtration.

The equivalence of categories now gives a natural transformation
of crossed complexes
\begin{equation} %\label{eta2}
\eta' \;:\; \Pi X_*  \otimes \Pi Y_*  \to \Pi(X_*\otimes Y_*).
\end{equation}It would be hard to construct this directly.
It is proved in \cite{bh:class91} that $\eta'$ is an isomorphism if
$X_*, Y_*$ are the skeletal filtrations of $CW$-complexes. The proof
uses the HHvKT, and the fact that $A \otimes -$ on crossed complexes
has a right adjoint and so preserves colimits. It is proved in
\cite{babr} that $\eta$ is an isomorphism if $X_*, Y_*$ are
cofibred, connected filtered spaces. This applies in particular to
the useful case of the filtration $B^\Box  C^* $ of the classifying
space of a crossed complex.

It turns out that the defining rules for the tensor product of
crossed complexes which follows from the above construction are
obtained as follows. We first define a bimorphism of crossed
complexes.

\begin{Def} \label{tensor} A \emph{bimorphism} $\theta : (A,B) \to C$
of crossed complexes is a family of maps $\theta : A_m \times B_n
\to C_{m+n}$ satisfying the following conditions, where $a \in
A_m, b \in B_n, a_1 \in A_1, b_1 \in B_1$ (temporarily using
additive notation throughout the definition):
\begin{enumerate}[(i)]
\item
  \begin{align*}
  \beta(\theta(a,b)) \;=\; \theta(\beta a,\beta b) ~
   \text{ for all } ~a \in A, b \in B\;.
  \end{align*}
\item
  \begin{align*}
  \theta(a,b^{b_1}) \;=\; \theta(a,b)^{\theta(\beta a,b_1)} ~
     \text{ if } ~m \ge 0, n \ge 2\;, \\
  \theta(a^{a_1},b) \;=\; \theta(a,b)^{\theta(a_1,\beta b)} ~
     \text{ if } ~m \ge 2, n \ge 0\;.
  \end{align*}
\item
  \begin{align*}
    \theta(a,b+b') & \;=\;
    \begin{cases}
      ~\theta(a,b) + \theta(a,b')
          &  \text{if } m=0,n \ge 1\
       \text{or } m \ge 1, n \ge 2\;, \\
      ~\theta(a,b)^{\theta(\beta a,b')} + \theta(a,b')
          &  \text{if } m \ge 1, n=1\;,
    \end{cases} \\
    \theta(a+a',b) & \;=\;
    \begin{cases}
      ~\theta(a,b) + \theta(a',b)
          &  \text{if } m \ge 1,n=0 \ \text{ or } m \ge 2, n \ge 1\;, \\
      ~\theta(a',b) + \theta(a,b)^{\theta(a',\beta b)}
          &  \text{if } m=1,n \ge 1\;.
    \end{cases}
  \end{align*}
\item
  \begin{align*}
  \delta_{m+n}(\theta(a,b)) & \;=\;
    \begin{cases}
    ~~\theta(\delta_m a,b) + (-)^m
     \theta(a,\delta_n b)
       &  \text{if } m \ge 2, n \ge 2\;,\\
    ~-\theta(a,\delta_n b) - \theta(\beta a,b) +
     \theta(\alpha a,b)^{\theta(a,\beta b)}
       &  \text{if } m=1,n \ge 2\;, \\
    ~(-)^{m+1} \theta(a,\beta b) + (-)^m
     \theta(a,\alpha b)^{\theta(\beta a,b)} + \theta(\delta_m a,b)
       &  \text{if } m \ge 2, n=1\;, \\
    ~-\theta(\beta a,b) - \theta(a,\alpha b) + \theta(\alpha a,b)
     +\theta(a, \beta b)
       &  \text{if } m=n=1\;.
    \end{cases}
  \end{align*}
\item
  \begin{align*}
  \delta_{m+n}(\theta(a,b)) & \;=\;
    \begin{cases}
    ~\theta(a,\delta_n b) & \text{if } m=0, n \ge 2\;, \\
    ~\theta(\delta_m a,b) & \text{if } m \ge 2, n=0\;.
    \end{cases}
  \end{align*}
\item
  \begin{align*}
  \alpha(\theta(a,b)) \,=\, \theta(a,\alpha b)
  \quad\text{and}\quad
   \beta(\theta(a,b)) \,=\, \theta(a,\beta b)
      & \quad\text{if}~ m=0,n=1\;, \\
  \alpha(\theta(a,b)) \,=\, \theta(\alpha a,b)
  \quad\text{and}\quad
   \beta(\theta(a,b)) \,=\, \theta(\beta a,b)
      & \quad\text{if}~ m=1,n=0\;.
  \end{align*}
\end{enumerate}

The \emph{tensor product} of crossed complexes $A,B$ is given by
the universal bimorphism $(A,B) \to A \otimes B$, $(a,b) \mapsto a
\otimes b$. The rules for the tensor product are obtained by
replacing  $\theta(a,b)$ by $a \otimes b$ in the above formulae.
\end{Def}

The conventions for these formulae for the  tensor product arise
from the derivation of the tensor product via the category of
cubical $\omega$-groupoids with connections, and the formulae are
forced by our conventions for the equivalence of the two
categories \cite{bh:algcub,bh:tens}.

The complexity of these formulae is directly related to the
complexities of the cell structure of the product $E^m \times E^n$
where the $n$-cell $E^n$ has cell structure $e^0$ if $n=0$,
$e^0_\pm \cup e^1$ if $n=1$, and $e^0 \cup e^{n-1} \cup e^n$ if $
n \ge 2$.

It is proved in \cite{bh:tens} that the bifunctor $- \otimes -$ is
symmetric and that if $a_0$ is a vertex of $A$ then the morphism
$B \to A \otimes B, \; b \to a_0 \otimes b$, is injective.

There is a standard groupoid model $\I$ of the unit interval, namely
the indiscrete groupoid on two objects $0,1$. This is easily
extended trivially to either a crossed complex or an \ogpd. So using
$\otimes$ we can define a `cylinder object' $\I \otimes -$ in these
categories and so a homotopy theory, \cite{br-go}.

\section{Free crossed complexes and free crossed resolutions}
\label{S:freecrossed}  Let $C$ be a crossed complex. A {\it free
basis} $B_*$ for $C$ consists of the following:

\noindent $B_0$ is set which we take to be $C_0$;

\noindent  $B_1$ is a graph with source and target maps $s,t:B_1
\to B_0$ and $C_1$ is the free groupoid on the graph $B_1$: that
is $B_1$ is a subgraph of $C_1$ and any graph morphism $B_1 \to G$
to a groupoid $G$ extends uniquely to a groupoid morphism $C_1 \to
G$;

\noindent $B_n$ is, for $n\ge 2$,  a totally disconnected subgraph
of $C_n$ with target map $t:B_n \to B_0$;  for $n=2$, $C_2$ is the
free crossed $C_1$-module on $B_2$ while for $n > 2$, $C_n$ is the
free $(\pi_1C)$-module on $B_n$.

It may be proved using the HHvKT that if $X_*$ is a $CW$-complex
with the skeletal filtration, then $\Pi X_* $ is the free crossed
complex on the characteristic maps of the cells of $X_*$. It is
proved in \cite{bh:class91} that the tensor product of free crossed
complexes is free.

A {\it free crossed resolution } $F_*$ of a groupoid $G$ is a free
crossed complex which is aspherical together with an isomorphism
$\phi: \pi_1(F_*) \to G$. Analogues of standard methods of
homological algebra show that free crossed resolutions of a group
are unique up to homotopy equivalence.

In order to apply this result to free crossed resolutions, we need
to replace free crossed resolutions by $CW$-complexes. A
fundamental result for this is the following, which goes back to
Whitehead \cite{wjhc:sht} and Wall \cite{wall}, and which is
discussed further by Baues in \cite[Chapter VI, \S 7]{baues1}:

\begin{thm}
Let $X_*$ be a $CW$-filtered space, and let $\phi: \Pi X_* \to C$ be
a homotopy equivalence to a free crossed complex with a preferred
free basis. Then there is a $CW$-filtered space $Y_*$, and an
isomorphism $\Pi Y_* \cong C$ of crossed complexes with preferred
basis, such that $\phi$ is realised by a homotopy equivalence $X_*
\to Y_*$.
\end{thm}

In fact, as pointed out by Baues, Wall states his result in terms
of chain complexes, but the crossed complex formulation seems more
natural, and avoids questions of realisability in dimension $2$,
which are unsolved for chain complexes.

\begin{cor} \label{cwmodel}
If $A$ is a free crossed resolution of a group $G$, then $A$ is
realised as free crossed complex with preferred basis by some
$CW$-filtered space $Y_*$.
\end{cor}
\begin{proof}
We only have to note that the group $G$ has a classifying
$CW$-space $BG$ whose fundamental crossed complex $\Pi( BG)$ is
homotopy equivalent to $A$.
\end{proof}

Baues also points out in \cite[p.657]{baues1} an extension of these
results which we can apply to the realisation of morphisms of free
crossed resolutions. A new proof of this extension is given by Faria
Martins in \cite{martins-bauesthm}, using methods of Ashley
\cite{ashley}.

\begin{prop} \label{cwmaps}
Let $X = K(G,1),\, Y = K(H,1)$ be $CW$-models of Eilenberg - Mac
Lane spaces and let $h : \Pi X_*  \to \Pi(Y_*)$ be a morphism of
their fundamental crossed complexes with the preferred bases given
by skeletal filtrations. Then  $h = \Pi(g)$  for some cellular $g: X
\to Y$.
\end{prop}
\begin{proof}
Certainly $h$ is homotopic to  $\Pi(f)$ for some $f : X \to Y$
since the set of pointed homotopy classes $X \to Y$ is bijective
with the morphisms of groups $A \to B$. The result follows from
\cite[p.657,(**)]{baues1} (`if $f$ is $\Pi$-realisable, then each
element in the homotopy class of $f$ is $\Pi$-realisable').
\end{proof}

These results are  exploited in \cite{emmathesis,brmopw} to
calculate free crossed resolutions of the fundamental groupoid of a
graph of groups.

An algorithmic approach to the calculation of free crossed
resolutions for groups is given in \cite{BRazak:LMS99}, by
constructing partial contracting homotopies for the universal cover
at the same time as constructing this universal cover inductively.
This has been implemented in GAP4 by Heyworth and Wensley
\cite{heywens}.

\section{Classifying spaces and the homotopy classification of maps}
\label{S:homclass} The formal relations of cubical sets and of
cubical $\omega$-groupoids with connections and the relation of Kan
cubical sets with topological spaces, allow the proof  of a homotopy
classification theorem:
\begin{thm}
If $K$ is a cubical set, and $G$ is an $\omega$-groupoid, then
there is a natural bijection of sets of homotopy classes
$$[|K|,|UG|] \cong [\rho(|K|_*), G],$$where on the left hand side we
work in the category of spaces, and on the right in
$\omega$-groupoids.
\end{thm}Here $|K|_*$ is the filtration by skeleta of
the geometric realisation of the cubical set.

We explained earlier how to define a cubical classifying space say
$B^\Box(C)$ of a crossed complex $C$ as $B^\Box(C)=  |U N^\Box C|=
|U\lambda C|$. The properties already stated now give the homotopy
classification theorem 4.1.9.

It is shown in \cite{bh:colimits} that for a  $CW$-complex $Y$ there
is a map $p:Y \to B^\Box\Pi Y_*$ whose homotopy fibre is
$n$-connected if $Y$ is connected and $\pi_i Y = 0$ for $2\le i \le
\ n - 1$. It follows that if also $X$  is a connected $CW$-complex
with $\dim X \le n$, then $p$ induces a bijection $$[X,Y] \to
[X,B\Pi Y_*].$$ So under these circumstances we get a bijection
\begin{equation}[X,Y] \to [\Pi X_*,\Pi Y_*]. \label{olum}\end{equation}
This result, due to Whitehead \cite{jhcw:CHII}, translates a
topological homotopy classification problem to an algebraic one.
We explain below how this result can  be translated to a result on
chain complexes with operators.

It is also possible to define a simplicial nerve $N^\Delta(C)$ of a
crossed complex $C$ by $$N^\Delta(C)_n= \crs(\Pi(\Delta^n),C).$$The
{\it simplicial classifying space} of $C$ is then defined using the
simplicial geometric realisation $$B^\Delta(C)= |N^\Delta(C)|.$$ The
properties of this simplicial classifying space are developed in
\cite{bh:class91}, and in particular an analogue of 4.1.9 is proved.

The simplicial nerve and an adjointness
$$\crs(\Pi(L),C) \cong \simp(L, N^\Delta(C))$$
are used in \cite{BGPT:Iindag,BGPT:IIK-theory} for an equivariant
homotopy theory of crossed complexes and their classifying spaces.
Important ingredients in this are notions of coherence and an
Eilenberg-Zilber type theorem for crossed complexes proved in Tonks'
Bangor thesis \cite{tonksthesis,tonks-ez}. See also \cite{rb-rs}.

Labesse in \cite{labesse} defines  a {\it crossed set}. In fact a
crossed set is exactly a crossed module (of groupoids) $\delta:C \to
X\rtimes G $ where $G$ is a group acting on the set $X$, and
$X\rtimes G$ is the associated actor groupoid; thus the simplicial
construction from a crossed set described by Larry Breen in
\cite{labesse} is exactly the simplicial nerve of the crossed
module, regarded as a crossed complex. Hence the cohomology with
coefficients in a crossed set used in \cite{labesse} is a special
case of cohomology with coefficients in a crossed complex, dealt
with in \cite{bh:class91}. (We are grateful to  Breen  for pointing
this out to us in 1999.)

\section{Relation with chain complexes with a groupoid of
operators} \label{S:chn} Chain complexes with a group of operators
are a well known tool in algebraic topology, where they arise
naturally as the chain complex $C_{*}\tilde{X}_*$ of cellular
chains of the universal cover $\tilde{X}_*$ of a reduced
$CW$-complex $X_*$. The group of operators here is the fundamental
group of the space $X$.

J.H.C. Whitehead in \cite{jhcw:CHII} gave an interesting relation
between his free crossed complexes (he called them `homotopy
systems') and such chain complexes. We refer later to his
important homotopy classification results in this area. Here we
explain the relation with the Fox free differential calculus
\cite{fox:1}.

Let $\mu : M \to P$ be a crossed module of groups, and let $G=
\Coker \mu$. Then there is an associated diagram
\begin{equation} \label{derived}
\xymatrix{ M \ar [r] ^\mu \ar [d]_{h_2} & P \ar [d] ^{h_1} \ar [r]
^\phi & G \ar [d] ^{h_0}\\
M ^\ab_{\phantom{0}} \ar [r] _{\partial _2} & D_\phi
^{\phantom{0}}\ar [r] _{\partial _1} & \Z[G] }\end{equation} in
which the second row consists of (right) $G$-modules and module
morphisms. Here $h_2$ is simply the Abelian isation map; $h_1: P \to
D_\phi$ is the universal $\phi$-derivation, that is it satisfies
$h_1(pq)=h_1(p)^{\phi q} + h_1(q)$, for all $p,q \in P$, and is
universal for this property; and $h_0$ is the usual derivation $g
\mapsto g-1$. Whitehead in his Lemma 7 of \cite{jhcw:CHII} gives
this diagram in the case $P$ is a free group, when he takes $D_\phi$
to be the free $G$-module on the same generators as the free
generators of $P$. Our formulation, which uses the derived module
due to Crowell \cite{crowell}, includes his case. It is remarkable
that diagram \eqref{derived} is a commutative diagram in which the
vertical maps are operator morphisms, and that the bottom row is
defined by this property. The proof in \cite{bh:chn} follows
essentially Whitehead's proof. The bottom row is exact: this follows
from results in \cite{crowell}, and is a reflection of a classical
fact on group cohomology, namely the relation between central
extensions and the Ext functor, see \cite{maclane:hom}.  In the case
the crossed module is the crossed module $\delta: C(\omega) \to
F(X)$ derived from a presentation of a group, then $C(\omega)^\ab$
is isomorphic to the free $G$-module on $R$, $D_\phi$ is the free
$G$-module on $X$, and it is immediate from the above that
$\partial_2$ is the usual derivative $(\partial r /\partial x)$ of
Fox's free differential calculus \cite{fox:1}. Thus Whitehead's
results anticipate those of Fox.

It is also proved in \cite{jhcw:CHII} that if the restriction $M
\to \mu(M)$ of $\mu$ has a section which is a morphism but not
necessarily a $P$-map, then $h_2$ maps $\Ker \mu $ isomorphically
to $\Ker \partial _2$. This allows calculation of the module of
identities among relations by using module methods, and this is
commonly exploited, see for example \cite{ellis-kholod} and the
references there.

Whitehead introduced the categories $\mathsf{CW}$ of reduced
$CW$-complexes, $\mathsf{HS}$ of homotopy systems, and
$\mathsf{FCC}$ of free chain complexes with a group of operators,
together with functors
\begin{equation*}
\mathsf{CW} \labto{\Pi} \mathsf{HS} \labto{C} \mathsf{FCC}.
\end{equation*}

In each of these categories he introduced notions of homotopy and
he proved that $C$ induces an equivalence of the homotopy category
of $\mathsf{HS}$ with a subcategory of the homotopy category of
$\mathsf{FCC}$. Further, $C\Pi X_*$ is isomorphic to the chain
complex $C_{*}\tilde{X}_*$ of cellular chains of the universal
cover of $X$, so that under these circumstances there is a
bijection of sets of homotopy classes
\begin{equation}
[\Pi X_*,\Pi Y_*] \rightarrow [C_{*} \tilde{X}_*, C_{*}
\tilde{Y}_*]. \label{piclass} \end{equation} This with the
bijection  \eqref{olum}  can be interpreted as an operator version
of the Hopf classification theorem. It is surprisingly little
known. It includes results of Olum \cite{olum:1953} published
later, and it enables quite useful calculations to be done easily,
such as the homotopy classification of maps from a surface to the
projective plane \cite{ellis:homclass}, and other cases. Thus we
see once again that this general theory leads to specific
calculations.

All these results are generalised in \cite{bh:chn} to the non free
case and to the non reduced case, which requires a groupoid of
operators,  thus giving functors
\begin{equation*} \FTop \labto{\Pi} \crs \labto{\nabla} \Chain
.\end{equation*} (The paper \cite{bh:chn} uses the notation $\Delta$
for this $\nabla$.) One utility of the generalisation to groupoids
is that the functor $\nabla$ then has a right adjoint, and so
preserves colimits. An example of this preservation is given in
\cite[Example 2.10]{bh:chn}. The construction of the right adjoint
$\Theta$ to $\nabla$ builds on a number of constructions used
earlier in homological algebra.

The definitions of the categories under consideration in order to
obtain a generalisation of the bijection \eqref{piclass} has to be
quite careful, since it works in the groupoid case, and not all
morphisms of the chain complex are realisable.

This analysis of the relations between these two categories is
used in \cite{bh:class91} to give an account of cohomology with
local coefficients.

It is also proved in \cite{bh:chn} that the functor $\nabla$
preserves tensor products, where the tensor in the category $\Chain$
is a generalisation to modules over groupoids of the usual tensor
for chain complexes of modules of groups. Since the tensor product
is described explicitly in dimensions $\le 2$ in \cite{bh:tens}, and
$(\nabla C)_n = C_n$ for $n \ge 3$,  this preservation yields  a
complete description of the tensor product of crossed complexes.

\section{Crossed complexes and simplicial groups and
groupoids}\label{S:simplicial}  The Moore complex $NG$ of a
simplicial group $G$ is not in general  a (reduced) crossed complex.
Let $D_nG$ be the subgroup of $G_n$ generated by degenerate
elements. Ashley showed in his thesis \cite{ashley} that $NG$ is a
crossed complex if and only if $(NG)_n \cap (DG)_n= \{1\}$ for all
$n \ge 1$.

Ehlers and Porter in \cite{ehlport1,ehlport2} show that there is a
functor $C$ from simplicial groupoids to crossed complexes in
which $C(G)_n$ is obtained from $N(G)_n$ by factoring out
$$(NG_n \cap D_n)d_{n+1}(NG_{n+1} \cap D_{n+1}),$$
where the Moore complex is defined so that its differential comes
from the last simplicial face operator.

This is one part of an investigation into the Moore complex of a
simplicial group, of which the most general investigation is by
Carrasco and Cegarra in \cite{carrasco&cegarra}.

An important observation in \cite{tp:topol} is that if $N
\triangleleft G$ is an inclusion of a normal simplicial subgroup of
a simplicial group, then the induced morphism on components $\pi_0
(N) \to  \pi_0(G)$ obtains the structure of crossed module. This is
directly analogous to the fact that if $F \to E \to B$ is a
fibration sequence then the  induced morphism of fundamental groups
$\pi_1(F,x) \to \pi_1(E,x)$ also obtains the structure of crossed
module. This last fact is relevant to algebraic $K$-theory, where
for a ring $R$ the homotopy fibration sequence is taken to be $F \to
B(GL(R)) \to B(GL(R))^+$.

\section{Other homotopy multiple groupoids}
\label{sec:otherhomgpds}

A natural question is whether there are other useful forms of higher
homotopy groupoids. It is because the geometry of convex sets is so
much more complicated in dimensions $> 1$ than in dimension $1$ that
new complications emerge for the theories of higher order group
theory and of higher homotopy groupoids. We have different
geometries for example those of disks, globes, simplices, cubes, as
shown in dimension 2 in the following diagram.
\begin{center}
\begin{overpic}[scale=0.7]{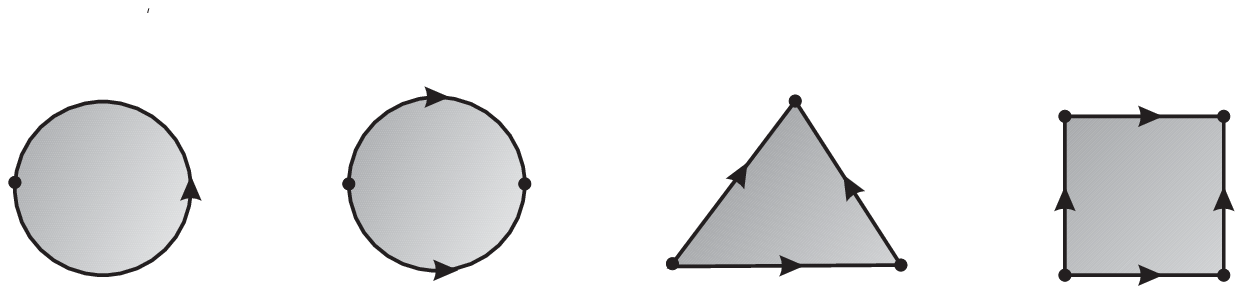}
\end{overpic}
\mbox{ }\\
\mbox{ }\\
\end{center}
The cellular decomposition for an $n$-disk is $D^n=e^0 \cup
e^{n-1}\cup e^n$, and that for globes is
$$G^n= e^0_\pm \cup e^1_\pm \cup \cdots \cup
e^{n-1}_\pm \cup e^n.$$ The higher dimensional group(oid) theory
reflecting the $n$-disks is that of crossed complexes, and that for
the $n$-globes is called {\it globular $\omega$-groupoids}.

A common notion of higher dimensional category is that of
$n$-category, which generalise the 2-categories studied in the late
1960s. A 2-category $\C$ is a category enriched in categories, in
the sense that each hom set $\C(x,y)$ is given the structure of
category, and there are appropriate axioms. This gives inductively
the notion of an $n$-category as a category enriched in
$(n-1)$-categories. This is called a `globular' approach to higher
categories. The notion of $n$-category for all $n$ was axiomatised
in \cite{bh1981a} and  called an $\infty$--category;   the
underlying geometry of a family of sets $S_n, n \geq 0$ with
operations $$D^\alpha_i:S_n \to S_{i}, E_i:S_{i} \to S_n,
\alpha=0,1; i=1,\ldots, n-1$$was there axiomatised. This was later
called a `globular set' \cite{street-globular}, and the term
$\omega$-category was used instead of the earlier $\infty$-category.
Difficulties of the globular approach are to define multiple
compositions, and also monoidal closed structures, although these
are clear in the cubical approach. A globular higher homotopy
groupoid of a filtered space has  been constructed in
\cite{brown:glob}, deduced from cubical results.

Although the proof of the HHvKT outlined earlier does seem to
require cubical methods, there is still a question of the place of
globular and simplicial   methods in this area. A simplicial
analogue of the equivalence of categories is given in
\cite{ashley,NanTie1}, using Dakin's notion of {\it simplicial
$T$-complex}, \cite{dakin}. However it is difficult to describe in
detail the notion of tensor product of such structures, or to
formulate a proof of the HHvKT theorem in that context. There is a
tendency to replace the term $T$-complex  from all this earlier work
such as \cite{bh:CR1,ashley} by {\it complicial set}, \cite{ver}.

It is easy to define a homotopy globular {\it set} $\rho^\bigcirc
X_* $ of a filtered space $X_*$ but it is not quite so clear how to
prove directly that the expected compositions are well defined.
However there is a natural graded map \begin{equation}i:
\rho^\bigcirc  X_* \to \rho X_* \label{globcub} \end{equation} and
applying the folding map of \cite{alagl:thesis,aa-b-s} analogously
to methods in \cite{bh:colimits} allows one to prove that $i$ of
\eqref{globcub} is injective. It follows that the compositions on
$\rho X_* $ are inherited by $\rho^\bigcirc  X_* $ to make the
latter a globular $\omega$-groupoid. The details are in
\cite{brown:glob}.

Loday in 1982 \cite{JLL2} defined the fundamental cat$^n$-group of
an $n$-cube of spaces (a cat$^n$-group may be defined as an $n$-fold
category internal to the category of groups), and showed that
cat$^n$-groups model all reduced weak homotopy $(n+1)$-types. Joint
work \cite{brlo:vkt} formulated and proved a HHvKT for the
cat$^n$-group functor from $n$-cubes of spaces. This allows new
local to global calculations of certain homotopy $n$-types
\cite{brown:adams}, and also an $n$-adic Hurewicz theorem,
\cite{brlo:hur}. This work obtains more powerful results than the
purely linear theory of crossed complexes. It yields a
group-theoretic description of the first non-vanishing homotopy
group of a certain $(n+1)$-ad of spaces, and so  several formulae
for the  homotopy and homology groups of specific spaces;
\cite{ellis-mikh} gives new applications. Porter in \cite{tp:topol}
gives an interpretation of Loday's results using methods of
simplicial groups. There is clearly a lot to do in this area. See
\cite{casas-etal} for relations of cat$^n$-groups with homological
algebra.

Recently some absolute homotopy 2-groupoids and double groupoids
have been defined, see \cite{BHKP:doublegroupoid} and the references
there, while \cite{brown-jan-double} applies generalised Galois
theory to give a new homotopy double groupoid of a map, generalising
previous work of \cite{bh1978}. It is significant that crossed
modules have been used in a differential topology situation by
Mackaay and Picken \cite{mack-pick1}. Reinterpretations of these
ideas in terms of double groupoids are started in \cite{br-gl1}.

It seems reasonable to suggest that in the most general case double
groupoids are still somewhat mysterious objects. The paper
\cite{andrus-natale} gives a kind of classification of them.

\section{Conclusion and questions}
\label{sec:conclusion}

$\bullet$ The emphasis on filtered spaces rather than the absolute
case is open to question.

$\bullet$  {\it Mirroring the geometry by the algebra} is crucial
for conjecturing and proving universal properties.

$\bullet$   {\it Thin elements} are crucial for modelling a concept
not so easy to define or handle algebraically, that of  commutative
cubes. See also \cite{hig-thin,steiner-thin}.

$\bullet$ The cubical methods summarised in section
\ref{sec:omegagroupoid} have also been applied in concurrency
theory, see for example \cite{gau-gou,fajetal}.

$\bullet$   {\it HHvKT theorems} give, when they apply,
 {exact information even in non commutative situations.} The
implications of this for homological algebra could be important.

$\bullet$  One construction inspired eventually by this work, the
{\it non Abelian  tensor product of groups}, has a bibliography of
90 papers since it was defined with Loday in \cite{brlo:vkt}.

$\bullet$ Globular methods do fit into this scheme. They have not so
far  yielded new calculations in homotopy theory, see
\cite{brown:glob}, but have  been applied to directed homotopy
theory, \cite{gau-gou}. Globular methods are the main tool in
approaches to weak category theory, see for example
\cite{leinster,street-globular}, although the potential of cubical
methods in that area is hinted at in \cite{steiner-thin}.

$\bullet$ For computations we really need strict structures
(although we do want to compute invariants of homotopy colimits).

$\bullet$ No work seems to have been done on Poincar\'e duality,
i.e. on finding special qualities of the fundamental crossed complex
of the skeletal filtration of a combinatorial manifold. However the
book by Sharko, \cite[Chapter VI]{sharko}, does use crossed
complexes for investigating Morse functions on a manifold.

$\bullet$ In homotopy theory, identifications in low dimensions
can affect high dimensional homotopy. So we need structure in a
range of dimensions to model homotopical identifications
algebraically. The idea of identifications in low dimensions is
reflected in the algebra by `induced constructions'.

$\bullet$ In this way we calculate some crossed modules modelling
homotopy 2-types, whereas the corresponding $k$-invariant is often
difficult to calculate.

$\bullet$ The use of crossed complexes in \v{C}ech theory is a
current project with Jim Glazebrook and Tim Porter.

$\bullet$  {\bf Question:} Are there applications of higher
homotopy groupoids  in other contexts where the fundamental
groupoid is currently used, such as algebraic geometry?

$\bullet$   {\bf Question:} There are uses of double groupoids in
differential geometry, for example in Poisson geometry, and in
2-dimensional holonomy \cite{br-icen:2-dhol}. Is there a non Abelian
De Rham theory, using an analogue of crossed complexes?

$\bullet$   {\bf Question:} Is there a truly non commutative
integration theory based on limits of multiple compositions of
elements of multiple groupoids?

\end{document}